\pdfoutput=1
\documentclass[a4paper]{article}
\usepackage[a4paper,top=3.5cm,bottom=2cm,left=2.5cm,right=2.5cm,marginparwidth=1.75cm]{geometry}

\usepackage{graphicx}
\usepackage{subcaption}
\usepackage[percent]{overpic}
\usepackage{amssymb}
\usepackage{amsmath, mathtools}
\usepackage[dvipsnames]{xcolor}
\usepackage{amsthm}
\usepackage{algorithm}
\usepackage[noEnd=true,endLComment=]{algpseudocodex}
\theoremstyle{remark}
\newtheorem*{remark}{\textbf{Remark}}
\usepackage[shortlabels]{enumitem}
\usepackage{nicefrac}
\usepackage[hidelinks]{hyperref}
\usepackage{multirow}
\usepackage[numbers, sort&compress]{natbib}
\usepackage[percent]{overpic}

\usepackage{booktabs}
\usepackage{makecell}

\newcommand{\jump}[1]{{\left[\!\left[#1\right]\!\right]}}
 
\newcommand{\ve}[1]{\text{$\boldsymbol{#1}$}} % vectors
\newcommand{\grad}{{\nabla}}

\renewcommand{\div}{\nabla\cdot}
\newcommand{\id}{\mathbb{I}}
\newcommand{\abs}[1]{\lvert#1\rvert}
\newcommand{\COMMENTBL}[1]{}
\newcommand{\convectiveterm}[2]{\left(#2\cdot \nabla\right)#1} % velocity, then convective velocity
\newcommand{\central}[1]{\{\!\{#1\}\!\}}

\ifpdf
  \DeclareGraphicsExtensions{.eps,.pdf,.png,.jpg}
\else
  \DeclareGraphicsExtensions{.eps}
\fi

\newcommand{\TheTitle}{A Discontinuous Galerkin Consistent Splitting Method for the Incompressible Navier--Stokes Equations}

\title{{\TheTitle}}
\author{Dominik T. Still \thanks{Institute for Computational Mechanics,
TUM School of Engineering and Design, Technical University of Munich
  (\texttt{dominik.still@tum.de}).} \footnotemark[3]
\and Natalia Nebulishvili \thanks{
Ruhr University Bochum - Faculty of Mathematics 
  (\texttt{natalia.nebulishvili@rub.de}, \texttt{richard.schussnig@rub.de}, \texttt{k.kormann@rub.de}, \texttt{martin.kronbichler@rub.de}).} \thanks{These authors contributed equally to this work.} 
\and Richard Schussnig\footnotemark[2]
\and Katharina Kormann\footnotemark[2]
\and Martin Kronbichler\footnotemark[2]
}

\begin{document}

\maketitle

\begin{abstract}
   This work presents the discontinuous Galerkin discretization of the consistent splitting scheme proposed by Liu [J. Liu, J. Comp. Phys., 228(19), 2009]. The method enforces the divergence-free constraint implicitly, removing velocity--pressure compatibility conditions and eliminating pressure boundary layers. Consistent boundary conditions are imposed, also for settings with open and traction boundaries. Hence, accuracy in time is no longer limited by a splitting error.
    
    The symmetric interior penalty Galerkin method is used for second spatial derivatives.
    The convective term is treated in a linearly implicit manner, which relaxes the CFL restriction of explicit schemes while avoiding the need to solve nonlinear systems required by fully implicit formulations. For improved mass conservation, Leray projection is combined with divergence and normal continuity penalty terms.
    By selecting appropriate fluxes for both the divergence of the velocity field and the divergence of the convective operator, the consistent pressure boundary condition can be shown to reduce to contributions arising solely from the acceleration and the viscous term for the $L^2$ discretization.
    
    Per time step, the decoupled nature of the scheme with respect to the velocity and pressure fields leads to a single pressure Poisson equation followed by a single vector-valued convection-diffusion-reaction equation.

    We verify optimal convergence rates of the method in both space and time and demonstrate compatibility with higher-order time integration schemes. A series of numerical experiments, including the two-dimensional flow around a cylinder benchmark and the three-dimensional Taylor--Green vortex problem, verify the applicability to practically relevant flow problems.

\end{abstract}

\noindent\textbf{Keywords.}
  time-splitting methods, consistent splitting, discontinuous Galerkin, fractional-step scheme, IMEX, incompressible Navier--Stokes \\

\section{Introduction}
\label{sec:introduction}

Incompressible fluid flows play a major role in science and engineering,
including countless applications of immense societal value. 
Computational fluid dynamics aims to provide approximate solutions with the help of computer simulation and has become
a fundamental discipline, immediately transferring to other areas of 
research and technology. Great progress has been made over the past decades, but even the
high-performance implementations of high-fidelity algorithms available today continue to 
challenge today's largest compute clusters. 
In a decades-long endeavor, accuracy and throughput of related numerical 
algorithms and their implementation have been continuously improved. 
However, further advancements are still urgently needed due to increasing demand in terms of accuracy, 
time to solution, numerical robustness, or conservation properties and physics conformity.
The innovations are often driven by challenging high-resolution applications in the energy, 
aviation or combustion sectors, or motivated by the requirements of many-query applications
such as uncertainty quantification, sensitivity analysis, model order reduction, inverse problem analysis and
machine learning.

The Navier--Stokes equations for transient incompressible flows are typically advanced in time by implicit
integrators due to the divergence-free or incompressibility constraint. 
Via the divergence-free condition, the pressure acts as a Lagrange multiplier enforcing
the continuity of mass. This structure calls for specific
choices of spaces for velocity and pressure approximation to fulfill the so-called inf-sup
(Ladyzhenskaya--Babu\v{s}ka--Brezzi) condition~\cite{Boffi2013, John2016}. 
Since the ground-breaking work by~\citet{TaylorHood1974}, who first used inf-sup stable
velocity--pressure pairs, 
the finite element method has been established as a viable alternative to the finite 
volume method. For inf-sup stability and 
suitable stabilization, see for instance~\cite{Brezzi1984, Hughes1986, CODINA2001, Pacheco2021_a, Schussnig2021_e}.

The standard time discretization approach, considering the mixed problem governing velocity and pressure fields as a large differential-algebraic system, leads to 
linear systems of saddle-point structure. Depending on the flow regime and problem parameters,
the pressure Schur complement can be approximated by pressure mass and Poisson 
operators~\cite{Cahouet1988}, while convective-dominant settings can be tackled via,
e.g., augmented-Lagrangian-based approaches~\cite{Benzi2006, Benzi2011, Heister2013, Lohmann2024augmented}.
However, these approaches can lead to an ill-conditioned velocity-velocity block, which in turn 
necessitates specialized multigrid methods, see, e.g.,~\cite{Shih2023}. 
Other---typically less Reynolds-robust alternatives---are the pressure 
convection-diffusion~\cite{Elman2014, Kay2003} 
and least-squares commutator~\cite{Elman2007, Elman2009} preconditioners,
which can perform well, especially for reasonable time step size and grid quality.
Some of these concepts have also been extended towards more complex physics, in the simplest case captured by nonlinear viscosity, used in, e.g., earth mantle convection~\cite{Rudi2015}, 
multiphase flows~\cite{Kronbichler2018}, or generalized Newtonian fluids~\cite{Schussnig2021_e}. Recent research has tried to re-formulate the partial differential equations underlying the saddle-point problem to simplify the linear algebra aspects~\cite{Creff2025}.

To (seemingly) avoid the complications associated with the saddle-point structure, a class of schemes, the so-called splitting, projection, or fractional-step schemes, circumvent solving a monolithic system at each time step by decoupling the governing equations for the velocity and pressure.
Note, however, that time splitting or projection schemes might still be subject to an inf-sup condition, limiting the number of admissible finite element pairings, see, e.g., \cite{CODINA2001, Fehn2017}.
Splitting the coupled problem into separate steps significantly simplifies the design of preconditioners, albeit
at the expense of more complex boundary conditions for the Poisson problem governing
the pressure. Relevant examples are incremental pressure-correction methods in rotational 
form as described in~\cite{Timmermans1996, Guermond2003rotational}, the dual splitting scheme 
by~\citet{Karniadakis1991}, see also \cite{Guermond2003a}, the consistent splitting 
schemes by~\citet{Guermond2003consistent} and~\citet{Liu2009} (see also the related method developed by~\citet{Huang2025stability}). Recent developments in this regard focus among other directions
on the non-trivial extensions towards complex rheological 
laws~\cite{Deteix2019, Pacheco2021_c, Schussnig2025b}, unconditional stability beyond 
second-order accurate schemes following ideas by~\citet{Huang2025stability}, 
as well as multiphase flows~\cite{Pacheco2022_a, Pacheco2025} and fluid--structure interaction~\cite{Schussnig2021_d}.

In this context, the present article focuses on spatial discretizations of higher-order accurate discontinuous Galerkin (DG)
discretizations of the consistent splitting scheme by~\citet{Liu2009} and thorough testing of the resulting framework. 
The motivation for choosing the consistent splitting over classical projection schemes is to avoid pressure boundary layers and achieve higher-order accuracy in time, while avoiding solving the monolithic system arising in mixed formulations. In contrast to the approach of~\citet{Liu2009}, which employs continuous finite elements, the DG discretization offers high-order accuracy combined with improved robustness in convection-dominated regimes~\cite{Hesthaven2007}.
The closest competitors are shortly commented on in what follows.

First, compared to the dual splitting scheme by~\citet{Karniadakis1991} and its robust DG discretization presented by~\citet{Krank2016,Fehn2017,Fehn2018b}, the present scheme employs consistent boundary conditions for the pressure Poisson equation and enables equal-order $C^0$-continuous interpolation~\cite{Johnston2004, Pacheco2021_a} according to numerical evidence. Furthermore, it is an important step towards unconditional stability for higher-order time stepping schemes, as the method can be extended adopting ideas by~\citet{Huang2025stability}. Second, \citet{Rosales2021} adopted the consistent splitting scheme and employed higher-order implicit-explicit Runge--Kutta schemes. However, they introduced the vorticity $\omega = \grad \times \ve{u}$ as an additional unknown to enforce the so-called electric boundary conditions naturally, a complication we can avoid due to the fact that we want to enforce a different set of boundary conditions, as will be discussed shortly. Third, \citet{Li2021} present a higher-order method based on the consistent splitting idea for the variable-density Navier--Stokes equations. They employ a DG method with bound-preserving limiter for the density evolution and continuous Galerkin schemes for the pressure and momentum equations combined with fully explicit and implicit-explicit Runge--Kutta time steppers. Compared to their work, we do not consider variable viscosity, but aim to improve robustness with respect to convective effects via DG discretization of the momentum equation.

The paper is organized as follows: Sec.~\ref{sec:splitting_scheme_idea} presents the fundamental idea behind consistent splitting schemes, leading to the temporal discretization presented in Sec.~\ref{sec:temporal_discretization}. The $L^2$-conforming DG discretization is discussed in Sec.~\ref{sec:spatial_discretization}. The properties of the devised numerical methods are critically investigated in Sec.~\ref{sec:numerical_results}, including spatial and temporal convergence studies in Sec.~\ref{sec:numerical_results_manufatured} on exact solutions, and the practically relevant benchmarks of flow around a cylinder~\cite{Schaefer1996}, see Sec.~\ref{sec:numerical_results_cylinder}, and the 3D Taylor--Green vortex problem~\cite{TaylorGreen1937, Wang2013}, see Sec.~\ref{sec:numerical_results_vortex}. Concluding remarks are given in Sec.~\ref{sec:conclusion}.

\section{Consistent Splitting Scheme}
\label{sec:splitting_scheme_idea}
\subsection{The Navier--Stokes Equations for Incompressible Flows}
We consider the Navier--Stokes equations governing the incompressible flow of a viscous fluid within the spatial domain $\Omega \subset \mathbb{R}^d$ with $d= 2 \text{ or } 3$. They consist of the momentum balance equation
\begin{equation} \label{eq:momentum_original}
   \frac{\partial \mathbf{u}}{\partial t} + (\mathbf{u} \cdot \nabla) \mathbf{u} - \nu \Delta \mathbf{u}  + \nabla p = \ve{f} \text{ in } \Omega \times (0,T] ,
\end{equation}
where $\nu$ is the kinematic viscosity, and the continuity equation
\begin{equation}
    \div \ve{u} = 0 \text{ in } \Omega \times (0,T] .
\end{equation}
The system is initialized with a divergence-free initial condition
\begin{equation}
    \ve{u}(\ve{x}, t=0) = \ve{u}_0(\ve{x}) \text{ in } \Omega.
\end{equation}
The boundary $\Gamma = \partial \Omega$ is split into non-overlapping Dirichlet and Neumann boundaries, denoted by $\Gamma^D$ and $\Gamma^N$, respectively, such that $\Gamma = \Gamma^D \cup \Gamma^N$ and $\Gamma^D \cap \Gamma^N = \emptyset$. On the Dirichlet and Neumann boundary segments the following conditions hold
\begin{align}
    \ve{u} &= \ve{g} &&\text{ on } \Gamma^D \times [0,T],
    \\
    (\nu \grad\ve{u} - p \id) \ve{n} &= \ve{h} &&\text{ on } \Gamma^N \times [0,T].
    \label{eq:traction_condition}
\end{align}
Here, $\ve{n}$ denotes the unit outward normal vector. 

\subsection{Consistent Splitting}
Standard implicit discretization of the Navier--Stokes problem~\eqref{eq:momentum_original}--\eqref{eq:traction_condition} adopting a mixed finite element method leads to a nonlinear system in saddle-point form. An attractive alternative is found in consistent splitting schemes, which split the coupled velocity-pressure system, resulting in a Poisson problem for the pressure and a convection-diffusion problem for the momentum equation. As proposed by~\citet{Liu2009}, the consistent splitting method replaces the incompressibility constraint by a pressure Poisson equation (PPE),
\begin{equation} \label{eq:ppe}
    - \Delta p = \div(\ve{u} \cdot \nabla \ve{u}) - \div \ve{f} \qquad \text{ in } \Omega \times [0,T] ,
\end{equation}
by taking the divergence of Equation~\eqref{eq:momentum_original} and using the divergence of the time derivative and the viscous term being zero under suitable regularity assumptions. 

To derive the boundary condition for the PPE on the Dirichlet portion $\Gamma^D$, Equation~\eqref{eq:momentum_original} is dotted with the normal vector $\ve{n}$, yielding
\begin{equation}
    \nabla p \cdot \ve{n} = 
        \left[
            -\frac{\partial \mathbf{g}}{\partial t} - (\mathbf{g} \cdot \nabla) \mathbf{u} + \ve{f}
        \right] \cdot \ve{n} 
        - \nu \nabla \times (\nabla \times \ve{u}) \cdot \ve{n} \qquad \text{ on } \Gamma^D \times [0,T].
\end{equation}
Here, we use the vector identity $\Delta \ve{u} \equiv \nabla(\div \ve{u}) - \nabla \times (\nabla \times \ve{u})$ to realize the viscous term in curl form as advantageous for splitting methods~\cite{Guermond2006}.

On the Neumann boundary, the pressure $g_p$ follows from the traction boundary condition~\eqref{eq:traction_condition} by projecting it onto the outward normal direction, see, e.g.~\cite{Pacheco2021_a}
\begin{align}
     g_p &\coloneqq p = \nu \nabla \ve{u} : \left( \ve{n} \otimes \ve{n} - \id \right) - \ve{h}  \cdot \ve{n}     &&\text{ on } \Gamma^N \times [0,T] .
    \label{eq:bc_traction_pressure}
\end{align}
Alternatively, the pressure $g_p$ may be applied directly on the boundary if known. Note that this choice also aligns with the most practically relevant scenarios, as frequently the pressure value is to be prescribed on open boundaries such as outflow. For a detailed discussion on this topic for this splitting scheme see, e.g.,~\cite{schussnig2022}.

\section{Discretization in Time}
\label{sec:temporal_discretization}
Within this work, backward differentiation formulae (BDF) are applied for time integration.
The time interval of interest $(0,T]$ is decomposed into time steps of size $\Delta t$ from $t^n$ to $t^{n+1}$. The scheme is written according to~\citet{Liu2009}, while the presentation of the individual steps is chosen similarly to~\cite{Pacheco2021_c}.

\subsection{Consistent Splitting Method}

As a first variant of the method, we consider the variant where the momentum equation is solved first with an extrapolated pressure term on the right-hand side, after which the pressure Poisson equation is used to compute the pressure implicitly with the updated velocity.
\paragraph{Momentum Step}
For the time derivative of the velocity in the momentum equation, a BDF scheme of order $J$ is selected. An extrapolation of order $J_c$ is employed for the pressure term on the right-hand side with $J_c = J$ or $J_c = J-1$ to enhance numerical stability---especially relevant for the higher-order case.
The velocity $\ve{u}^{n+1}$ at time $t^{n+1}$ is computed by solving
\begin{equation} \label{eq:momentum_equation}
    \frac{\gamma_0}{\Delta t}
				\ve{u}^{n+1}
				+
				\convectiveterm{\ve{u}^{n+1}}{ \ve{u}^{\star} }
				- \nu
				\div \left( \,  \grad \ve{u}^{n+1}\right)
				=
				\ve{f}^{n+1} % \left(t_{n+1}\right)
				+
				\sum_{i=1}^{J}
				\frac{\alpha_i}{\Delta t}
				\ve{u}^{n+1-i}
				-
				\sum_{i=1}^{J_c} \Bar{\beta_i} \nabla p^{n+1-i},
\end{equation}%
where $\gamma_0$, $\alpha_i$ and $\beta_i$, $i=1, \dots, J$ and $\bar{\beta_i}$, $i = 1, \dots, J_c$ are suitable time integration and extrapolation constants (see, e.g.,~\citet{Hairer1993}). The extrapolated velocity is defined as $\ve{u}^{\star} \coloneqq \sum_{i=1}^{J} \beta_i \ve{u}^{n+1-i}$. This linearly implicit form of the convective term relaxes the strict CFL conditions given by treating convection fully explicitly, while also avoiding the need to solve a nonlinear system of equations as resulting from fully implicit time integration, yet all options produce valid schemes, see also~\cite{Pacheco2021_c} and \cite{Ascher1995} for more general discussions on implicit-explicit schemes. 

\paragraph{Pressure Step} With the velocity $\ve{u}^{n+1}$, the PPE is solved for the pressure $p^{n+1}$ via
\begin{align} \label{eq:ppe_time_discretized}
 - \Delta p^{n+1} &=
		 \div \left[\convectiveterm{\ve{u}^{n+1}}{\ve{u}^{n+1}} \right]	
		- \nabla \cdot \ve{f}^{n+1},
\end{align}%
with the consistent boundary condition 
\begin{equation} \label{eq:ppe_dirichlet}
	\begin{aligned}
   \nabla p^{n+1} \cdot \ve{n} =  & 
				 \left[ -\convectiveterm{\ve{u}^{n+1}}{\ve{g}^{n+1}} \right] \cdot \ve{n}
                 -  \left[
					 \nu \nabla \times \left( \nabla \times \ve{u}^{n+1} \right) \right] \cdot \ve{n}  \\
				& - \partial_t \ve{g}^{n+1} \cdot \ve{n} + \ve{f}^{n+1} % \left(t_{n+1}\right)  
                 \cdot \ve{n}
				\quad
				\text{on }\Gamma_\mathrm{D}.
    \end{aligned}
\end{equation}
If the time derivative of the boundary condition $\partial_t \ve{g}^{n+1}$ is not known analytically, a BDF scheme of order $J$ is employed for approximation~\cite{Fehn2017}.
On the Neumann boundary $\Gamma_\mathrm{N}$, we can directly apply the decomposition from Equation~\eqref{eq:bc_traction_pressure} to get
\begin{align}
     g_p^{n+1} = \nu \nabla \ve{u}^{n+1} : \left( \ve{n} \otimes \ve{n} - \id \right)  - \ve{h}^{n+1}  \cdot \ve{n}     &&\text{ on } \Gamma^N,
    \label{eq:bc_traction_pressure_time_discrete}
\end{align}
or apply $g_p^{n+1}$ directly if known.

\subsection{Leray Projection}
The divergence-free condition is only implicitly enforced when adopting the consistent splitting scheme by~\citet{Liu2009}, which is rooted in the PPE replacing the continuity equation to begin with. To additionally suppress the divergence in the velocity field, we employ a Helmholtz--Leray decomposition, 
\begin{equation}
    \ve{u} = \hat{\ve{u}} - \nabla \phi
    ,
    \label{eq:helmholtz_leray_decomposition}
\end{equation} 
see \citet{Liu2009, Pacheco2021_c} for a more detailed discussion.
In this presentation, $\hat{\ve{u}}$ is solenoidal (divergence-free at every point), while $\nabla \phi$ is irrotational. This splitting can be used to derive the so-called Leray projection by taking the divergence of Equation~\eqref{eq:helmholtz_leray_decomposition}
\begin{align}
    \nabla \cdot \ve{u} &= \nabla \cdot \hat{\ve{u}} - \nabla \cdot \nabla \phi
    \nonumber
    \\
    \Leftrightarrow \quad
     \nabla \cdot \ve{u} &=  - \nabla \cdot \nabla \phi
    \label{eq:leray_projection}
    .
\end{align}
Equation~\eqref{eq:leray_projection} is completed with homogeneous Dirichlet and Neumann boundary conditions
\begin{equation}
    \ve{n}\cdot\nabla \phi = 0
    \quad \text{on}~\Gamma_\mathrm{D}
    ,
    \qquad \mathrm{and} \qquad
    \phi = 0
    \quad \text{on}~\Gamma_\mathrm{N}.
    \label{eq:leray_boundary_conditions}
\end{equation}

The velocity $\hat{\ve{u}} = \ve{u} + \nabla \phi$ is then completely divergence-free. 
Note, however, that $\hat{\ve{u}}$ does not fulfill the tangential boundary conditions of the original problem. This can easily be shown by forming the dot product of the Helmholtz--Leray decomposition~\eqref{eq:helmholtz_leray_decomposition} with a tangential unit vector $\ve{\tau}$
which yields
\begin{equation*}
    \ve{\tau} \cdot \hat{\ve{u}}
    =
    \ve{\tau} \cdot \ve{u}
    +
    \ve{\tau} \cdot \nabla \phi
    .
\end{equation*}
Simultaneously, we obtain the correct normal boundary condition of $\ve{u}$
\begin{equation*}
    \ve{n} \cdot \hat{\ve{u}}
    =
    \ve{n} \cdot \ve{u}
    +
    \ve{n} \cdot \nabla \phi
    =
    \ve{n} \cdot \ve{u}
\end{equation*}
due to the boundary condition for the normal derivative of $\phi$ on $\Gamma_\mathrm{D}$, see Equation~\eqref{eq:leray_boundary_conditions}. This dilemma forces us to choose between a) a not strictly divergence-free velocity field with correct boundary conditions obtained via the standard form we started with given in Equations~\eqref{eq:momentum_equation}--\eqref{eq:ppe_time_discretized}, or b) a solenoidal velocity field $\hat{\ve{u}}$ which does not fulfill the tangential boundary conditions as desired. Following~\citet{Liu2009,Pacheco2021_c} among others, we apply the Leray projection only for the velocities of old time steps, which might be compared to eliminating the projection step in the context of classical projection schemes, see, e.g., \citet{Guermond2006} Equations (3.17)--(3.18). The divergence-free velocity field $\hat{\ve{u}}$ is employed in the acceleration
term on the right-hand side of Equation~\eqref{eq:momentum_equation}, which can be rewritten as

\begin{align} \label{eq:momentum_equation_leray}
    \frac{\gamma_0}{\Delta t}
				\ve{u}^{n+1}
				&+
				\convectiveterm{\ve{u}^{n+1}}{ \ve{u}^{\star} }
				- \nu
				\div \left( \,  \grad \ve{u}^{n+1}\right) 
				\nonumber
				\\
				&=
				\ve{f}^{n+1}
				+
				\sum_{i=1}^{J}
				\frac{\alpha_i}{\Delta t}
				\hat{\ve{u}}^{n+1-i}
				-
				\sum_{i=1}^{J_c} \bar{\beta_i}
				\nabla p^{n+1-i} 
				\nonumber
				\\
				&=
				\ve{f}^{n+1}
				+
				\sum_{i=1}^{J}
				\frac{\alpha_i}{\Delta t}
				( \ve{u}^{n+1-i} + \nabla \phi^{n+1-i})
				-
				\sum_{i=1}^{J_c} \bar{\beta_i}
                \nabla p^{n+1-i}.
\end{align}
The momentum equation including the Leray correction is solved first. With the new velocity, the PPE~\eqref{eq:ppe_time_discretized} yields the pressure $p^{n+1}$. The Leray variable $\phi$ is determined via
\begin{align} \label{eq:leray_time}
     -\Delta \phi^{n+1} = \nabla \cdot \ve{u}^{n+1},
\end{align}
with boundary conditions given in Equation~\eqref{eq:leray_boundary_conditions}. This leads to an additional Poisson equation to be solved per time step, where $\phi^{n+1}$ is used in the next time steps.

\subsection{Scheme with Modified Pressure Variable}
The consistent splitting scheme with Leray correction consists of Equations~\eqref{eq:momentum_equation}--\eqref{eq:ppe_time_discretized} and the additional PPE~\eqref{eq:leray_time}. This increased cost compared to the standard scheme can be eliminated by defining a modified pressure variable $\hat{p}$. Therefore, we start with the momentum Equation~\eqref{eq:momentum_equation_leray} with the modified acceleration term
\begin{align} \label{eq:momentum_equation_leray_modifed_pressure}
    \frac{\gamma_0}{\Delta t}
				\ve{u}^{n+1}
				&+
				\convectiveterm{\ve{u}^{n+1}}{ \ve{u}^{\star} }
				- \nu
				\div \left( \,  \grad \ve{u}^{n+1}\right) 
                \nonumber
				\\
				&=
				\ve{f}^{n+1}
				+
				\sum_{i=1}^{J}
				\frac{\alpha_i}{\Delta t}
				  \ve{u}^{n+1-i}
                +
				\sum_{i=1}^{J}
				\frac{\alpha_i}{\Delta t}
                \nabla \phi^{n+1-i}
				-
				\sum_{i=1}^{J_c} \bar{\beta_i}
                \nabla p^{n+1-i}
                \nonumber
				\\
				&=
				\ve{f}^{n+1}
				+
				\sum_{i=1}^{J}
				\frac{\alpha_i}{\Delta t}
				  \ve{u}^{n+1-i}
				-
                \nabla  \hat{p}^{n+1}.
\end{align}
Instead of computing $\phi$ separately, the modified pressure $\hat{p}^{n+1}$ is defined as
\begin{equation} \label{eq:modified_pressure}
    \hat{p}^{n+1} = \sum_{i=1}^{J_c} \bar{\beta_i} p^{n+1-i}
                  - \sum_{i=1}^{J}
				\frac{\alpha_i}{\Delta t} \phi^{n+1-i}.
\end{equation}
Note that the momentum Equation~\eqref{eq:momentum_equation_leray_modifed_pressure} is unchanged compared to Equation~\eqref{eq:momentum_equation_leray}. Instead, the Leray correction term and the extrapolated pressure term are combined into one variable.

To advance one time step, the modified pressure $\hat{p}^{n+1}$ is first computed by solving a modified pressure Poisson equation, which is derived by linearly combining (instantaneous) PPEs~\eqref{eq:ppe_time_discretized} at previous time steps and Equation~\eqref{eq:leray_time}. The modified PPE reads 
\begin{align} \label{eq:modified_ppe}
- \Delta \hat{p}^{n+1} &= -	\sum_{i=1}^{J_c} \bar{\beta_i}
				\Delta p^{n+1-i}  
				+ \sum_{i=1}^{J}
				\frac{\alpha_i}{\Delta t} 
				\Delta \phi^{n+1-i}	
				\nonumber
				\\
				&=
				\sum_{i=1}^{J_c} \bar{\beta_i} \div (
		\ve{u}^{n+1-i} \cdot \nabla \ve{u}^{n+1-i}
		-  \ve{f}^{n+1-i})
        - \sum_{i = 1}^{J} \frac{\alpha_i}{\Delta t} \div \ve{u}^{n+1-i},
\end{align}
where all involved terms are known from previous time steps.

A single time step thus comprises the following two steps: 
\begin{enumerate}
    \item solve the modified PPE~\eqref{eq:modified_ppe} for the modified pressure $\hat{p}^{n+1}$,
    \item solve the Leray corrected momentum Equation~\eqref{eq:momentum_equation_leray_modifed_pressure} for the new velocity $\ve{u}^{n+1}$, which is equivalent to Equation~\eqref{eq:momentum_equation_leray}.
\end{enumerate}
The sequence of steps is reversed from the one outlined above in the derivation, but still keeps the PPE implicit for each time step.
The modified pressure variable $\hat{p}^{n+1}$ is a consistent pressure approximation (cf. \cite{Liu2009, Pacheco2021_c}), and can hence be used directly. This eliminates the computational overhead otherwise introduced by the Leray projection: Instead of solving two Poisson problems---one for the pressure $p^{n+1}$ and one for the Leray variable $\phi^{n+1}$ after computing the velocity $\ve{u}^{n+1}$---the modified pressure $\hat{p}^{n+1}$ is computed before solving the momentum step, combining the contributions of the two Poisson problems. 

The term $\sum_i^{J_c} \bar{\beta_i}\ve{f}^{n+1-i} \simeq \ve{f}^{n+1}$ in Equation~\eqref{eq:modified_ppe} extrapolates the forcing term, such that the modified PPE can finally be written as
\begin{equation} \label{eq:ppe_time_discretized_leray}
    - \Delta \hat{p}^{n+1} =
				\sum_{i=1}^{J_c} \bar{\beta_i} \div (
		\ve{u}^{n+1-i} \cdot \nabla \ve{u}^{n+1-i}) - \nabla \cdot \ve{f}^{n+1}
        - \sum_{i = 1}^{J} \frac{\alpha_i}{\Delta t} \div \ve{u}^{n+1-i}.	
\end{equation}
In the consistent boundary condition the convective term is again extrapolated with order $J_c$, while for the viscous term an extrapolation order $J_p$ can be chosen independently. Again suitable candidates are $J_p = J$ or $J_p = J - 1$ (cf.~\cite{Liu2010,Karniadakis1999}). The modified PPE is then subject to the boundary conditions 
\begin{equation} \label{eq:ppe_dirichlet_time_discretized}
	\begin{aligned}
   \nabla \hat{p}^{n+1} \cdot \ve{n} =  &  \sum_{i=1}^{J_c} \Bar{\beta_i}
				 \left[ -\convectiveterm{\ve{u}^{n+1-i}}{\ve{g}^{n+1-i}} \right] 
                 + \sum_{i=1}^{J_p} \hat{\beta_i} \left[
					- \nu \nabla \times \left( \nabla \times \ve{u}^{n+1-i} \right) \right] \cdot \ve{n}  \\
				& - \partial_t \ve{g}^{n+1} \cdot \ve{n} + \ve{f}^{n+1} % \left(t_{n+1}\right)  
                 \cdot \ve{n}
				\quad
				\text{on }\Gamma_\mathrm{D},
    \end{aligned}
\end{equation}
with suitable extrapolation coefficients $\hat{\beta_i}$, $i=1,\dots,J_p$ (see, e.g., \citet{Hairer1993}). The effects of choosing different orders for $J$, $J_c$ and $J_p$ are later shown for a numerical example in Tab.~\ref{tab:taylor_green_bdf}. 

Note that due to the boundary conditions of $\phi$ given in Equation~\eqref{eq:leray_boundary_conditions}, Equation~\eqref{eq:ppe_dirichlet_time_discretized} is not affected and remains unaltered.
On the Neumann boundary the pressure $\hat{g}_p^{n+1}$ is given by
\begin{align}
     \hat{g}_p^{n+1} = \sum_{i=1}^{J_h} \Tilde{\beta_i} \left( \nu \nabla \ve{u}^{n+1-i} : \left( \ve{n} \otimes \ve{n} - \id \right) \right) - \ve{h}^{n+1}  \cdot \ve{n}     &&\text{ on } \Gamma^N.
    \label{eq:bc_traction_pressure_time_discrete_leray} 
\end{align}
with respective coefficients $\Tilde{\beta_i}$, $i=1,\dots,J_h$, yet again unaltered by Equation~\eqref{eq:leray_boundary_conditions}.

This approach follows the pressure-approximation form of~\citet{Liu2010}, who showed that the method achieves the convergence order $\mathcal{O}(\Delta t^J)$ for a BDF scheme of order $J$.
The only difference with respect to~\cite{Liu2010} is that the convective term in the modified PPE and the boundary condition is treated with order $J_c$ instead of $J$, giving more flexibility.

\subsection{Iterative Scheme}
In order to increase the accuracy and robustness of the method in case of strong nonlinearities, the scheme can be slightly adapted, iterating between the PPE and the momentum equation until convergence is reached. The iteration is initialized by solving Equations~\eqref{eq:ppe_time_discretized_leray} and~\eqref{eq:momentum_equation_leray_modifed_pressure} to obtain the starting value for the velocity $\Tilde{\ve{u}}$. One iteration consists of solving the pressure step for $\Tilde{p}$,
\begin{equation}
    \label{eq:pressure_iteration}
    - \Delta \Tilde{p} = \nabla \cdot \left(  (\Tilde{\ve{u}} \cdot \nabla) \Tilde{\ve{u}}  \right) - \nabla \cdot  \ve{f}^{n+1} - \sum_{i=1}^{J} \frac{\alpha_i}{\Delta t} \nabla \cdot \ve{u}^{n+1-i},
\end{equation}
with boundary conditions adopted from Equations~\eqref{eq:ppe_dirichlet}--\eqref{eq:bc_traction_pressure_time_discrete} and the momentum step, setting the convective velocity equal to the last velocity iterate, $\ve{u}^{\ast} = \Tilde{\ve{u}}$, and solving
\begin{align}
 \frac{\gamma_0}{\Delta t}
				\Tilde{\ve{u}}
				%&
                +
				\convectiveterm{\Tilde{\ve{u}}}{ \ve{u}^{\ast}}
				- \nu
				\div \left( \,  \grad \Tilde{\ve{u}}\right) 
				%\nonumber
				%\\
				%&
                =
				\ve{f}^{n+1}
				+
				\sum_{i=1}^{J}
				\frac{\alpha_i}{\Delta t}
				\ve{u}^{n+1-i}
				- \nabla \Tilde{p}
                .
                \label{eq:momentum_iteration}
\end{align}
Once convergence is reached, pressure and velocity are set to $p^{n+1} = \Tilde{p}$ and $\ve{u}^{n+1} = \Tilde{\ve{u}}$.

Note that the increased robustness and accuracy come at a significant price to solve several systems until convergence. Each iteration of this scheme requires solving the momentum step and the modified PPE. Naturally, whether this added cost per time step amortizes depends on the target application at hand, but the overwhelming majority of popular time splitting or projection methods are typically iteration-free for performance reasons. We therefore present the iterated variant primarily to gain further insight into the methods via a comparison in Section~\ref{sec:numerical_results_vortex}, and not as a particularly strong candidate due to its increased cost per time step. 
Compared to standard mixed solution approaches, though, saddle point systems are avoided in any case.

\subsection{Overview} \label{chap:variants}
Table~\ref{tab:overview_schemes} summarizes the consistent splitting schemes by specifying the sequence of momentum and pressure equations, including the associated boundary conditions, as applied within a single time step.

\begin{table}[h]
\centering
\caption{Versions of the consistent splitting scheme to advance one step in time: respective momentum and pressure steps, and related boundary conditions (BCs).}
\label{tab:overview_schemes}
\begin{tabular}{ll}
\toprule
Scheme & Steps \\
\midrule

Consistent Splitting &
\makecell[l]{
1. Momentum equation~\eqref{eq:momentum_equation}\\
2. PPE~\eqref{eq:ppe_time_discretized} 
%\\
%\quad 
with BCs: \eqref{eq:ppe_dirichlet}--\eqref{eq:bc_traction_pressure_time_discrete}
} \\
\midrule

Leray Projection &
\makecell[l]{
1. Corrected momentum equation~\eqref{eq:momentum_equation_leray}\\
2. PPE~\eqref{eq:ppe_time_discretized} 
%\\
%\quad
with BCs: \eqref{eq:ppe_dirichlet}--\eqref{eq:bc_traction_pressure_time_discrete} \\
3. Leray projection~\eqref{eq:leray_time}
} \\
\midrule

Modified pressure &
\makecell[l]{
1. Modified PPE~\eqref{eq:ppe_time_discretized_leray} 
%\\
%\quad 
with BCs: \eqref{eq:ppe_dirichlet_time_discretized}--\eqref{eq:bc_traction_pressure_time_discrete_leray} \\
2. Corrected momentum equation~\eqref{eq:momentum_equation_leray_modifed_pressure}
} \\
\midrule

Iterative scheme &
\makecell[l]{
1. Modified PPE~\eqref{eq:ppe_time_discretized_leray} 
%\\
%\quad 
with BCs: \eqref{eq:ppe_dirichlet_time_discretized}--\eqref{eq:bc_traction_pressure_time_discrete_leray} \\
2. Corrected momentum equation~\eqref{eq:momentum_equation_leray_modifed_pressure} \\
3. Iteration:\\
\quad\quad a) Modified PPE~\eqref{eq:pressure_iteration} 
%\\
%\quad\quad\quad 
with BCs: \eqref{eq:ppe_dirichlet}--\eqref{eq:bc_traction_pressure_time_discrete}\\
\quad\quad b) Corrected momentum equation~\eqref{eq:momentum_iteration}
} \\

\bottomrule
\end{tabular}
\end{table}

\subsection{Similarity to the Higher-Order Dual Splitting with Periodic Boundary Conditions}
In the following, similarities between the scheme with modified pressure variable and the higher-order dual splitting by \citet{Karniadakis1991} are shown. The dual splitting scheme is given by
\begin{align}
  \frac{\gamma_0 \mathbf{u}^* - \sum_{i = 0}^{J-1} (\alpha_i \mathbf{u}^{n-i})}{\Delta t} &= -\sum_{i = 0}^{J-1} \beta_i (\mathbf{u}^{n-i} \cdot \nabla) \mathbf{u}^{n-i} +  \mathbf{f}^{n+1}, \label{eq:dual_splitting_1}\\
   -\Delta p^{n+1} &= -\frac{\gamma_0}{\Delta t} \nabla \cdot \mathbf{u}^*, \label{eq:dual_splitting_pressure_poisson}\\ 
  \mathbf{u}^{**} &= \mathbf{u}^* - \frac{\Delta t}{\gamma_0} \nabla p^{n+1},  \label{eq:dual_splitting_3}\\
 \frac{\gamma_0}{\Delta t} \mathbf{u}^{n+1} - \nu \Delta \mathbf{u}^{n+1}  &= \frac{\gamma_0}{\Delta t} \mathbf{u}^{**} .\label{eq:dual_splitting_helmholtz}
\end{align}

Inserting Equation~\eqref{eq:dual_splitting_1} into~\eqref{eq:dual_splitting_pressure_poisson} leads to the PPE with Leray correction given by Equation~\eqref{eq:ppe_time_discretized_leray}. Further, if Equations~\eqref{eq:dual_splitting_3} and~\eqref{eq:dual_splitting_1} are inserted into Equation~\eqref{eq:dual_splitting_helmholtz}, it is similar to the momentum Equation~\eqref{eq:momentum_equation}. The difference is the treatment of the convective term, which is explicit in the dual splitting scheme and linearly implicit in the consistent splitting method presented here. Note, however, that in the dual splitting scheme it could also be considered linearly implicit, see the variant presented in~\cite{Liosi2025} and related work in~\cite{Guermond2003a, Schussnig2025b, Guermond2006, Sherwin2003}.

For an explicit treatment of the convective operator within the proposed consistent splitting method, the two schemes differ in the boundary conditions only, see Sec.~\ref{sec:Appendix}. Under periodic boundary conditions, the two schemes are identical, assuming that only the temporal discretization is considered, i.e., when the problem remains continuous in space.
For the proposed $L^2$-conforming spatial discretization, the right-hand side term of Equation~\eqref{eq:dual_splitting_1} can be treated with classical flux boundary conditions, whereas the right-hand side of the PPE incorporates the convective term tested against the gradient of the test function, see Equation~\eqref{eq:rhs_ppe_convective}.

This equivalence raises doubts about the requirements regarding the choice of function spaces for velocity and pressure in the consistent splitting method presented here. The dual splitting method of \citet{Karniadakis1991} is known to be subject to an inf-sup condition, while the consistent splitting scheme (without Leray projection) is reportedly not subject to such a condition~\cite{Liu2009,Liu2010,Liu2009c}. Adopting the Leray projection required to obtain the optimal convergence rates for the present formulation, however, might re-introduce the inf-sup condition in some sense. Hence, the potential accuracy loss using equal-order interpolation for viscous flows due to the curl-curl term is not the only reason for choosing $k_u = k_p + 1$.

Finally, we remark that the proposed consistent splitting scheme is also advantageous for $H^1$-conforming or $H(\text{div})$-conforming discretizations, since it avoids the solution of a mass matrix system for step~\eqref{eq:dual_splitting_1} present in the dual splitting method.
A comparison of the present splitting scheme and the dual splitting scheme by \citet{Karniadakis1991} is provided in Sec.~\ref{sec:numerical_results_vortex}.

\section{Discretization in Space}
\label{sec:spatial_discretization}
This section presents the discretization with an $L^2$-conforming discontinuous Galerkin method for the modified pressure scheme. The other schemes discussed in Section~\ref{sec:temporal_discretization} follow analogously, but are omitted here for brevity. To this end, the spatial domain $\Omega$ is approximated by $N_\mathrm{el}$ non-overlapping elements $\Omega_e$, $e=1,\dots,N_\mathrm{el}$, such that $\Omega \approx \Omega_h \coloneqq \bigcup_{e=1}^{N_\mathrm{el}} \Omega_e$. 
The spaces of the test and ansatz functions are defined as
{
    \footnotesize
    \begin{equation} \label{eq:function_spaces}
        \begin{aligned}
            \mathcal{V}_h^u &= 
            \biggl\{\, \ve{u}_h \in [L^2(\Omega_h)]^d :\,
            \ve{u}_h\bigl(\ve{x}(\ve{\xi})\bigr)\big|_{\Omega_e}
            = \widetilde{\ve{u}}_h^{\,e}(\ve{\xi})\big|_{\widetilde{\Omega}_e}
            \in V_{h,e}^u
            = \bigl[\mathbb{P}_{k_u}(\widetilde{\Omega}_e)\bigr]^d,\,
            \forall e=1,\dots,N_{\mathrm{el}} \,\biggr\}, \\
            \mathcal{V}_h^p &= 
            \biggl\{\, p_h \in L^2(\Omega_h) :\, 
            p_h\bigl(\ve{x}(\ve{\xi})\bigr)\big|_{\Omega_e}
            = \widetilde{p}_h^{\,e}(\ve{\xi})\big|_{\widetilde{\Omega}_e}
            \in V_{h,e}^p
            = \mathbb{P}_{k_p}(\widetilde{\Omega}_e),\,
            \forall e=1,\dots,N_{\mathrm{el}} \,\biggr\}
            ,
        \end{aligned}
    \end{equation}
}%
where $\ve{x}$ is the spatial coordinate in each element $\Omega_e$, connected to the reference element $\widetilde{\Omega}_e$ with reference coordinate $\ve{\xi}$ via an element-wise isoparametric mapping. On the reference element, the polynomial space $\mathbb{P}_{k}$ of degree $k$ is constructed via Lagrange polynomials. Within this work, inf-sup stable velocity-pressure pairings are considered, i.e., $k_u = k_p + 1$.
The functions $\ve u_h \in \mathcal{V}_h^u$ and $p_h \in \mathcal{V}_h^p$ approximate the velocity and the pressure, while the test functions from the velocity and pressure space are denoted by $\ve{v}_h\in \mathcal{V}_h^u$ and $q_h \in \mathcal{V}_h^p$, respectively. Note, however, that for simplicity of notation, the discrete test and trial functions (unless stated otherwise) are written as $\ve{u}$, $\ve{v}$, $p$ and $q$. We further define the jump and average operators along inner faces as
\begin{gather*}
	\jump{(\cdot)} \coloneqq (\cdot)_- - (\cdot)_+
    \qquad \text{and} \qquad
	\central{(\cdot)} \coloneqq \nicefrac{1}{2}\left((\cdot)_- + (\cdot)_+\right)
    ,
\end{gather*}
where $(\cdot)$ might be scalar or vector-valued. Here, $(\cdot)_-$ denotes quantities on the
current element, while $(\cdot)_+$ denotes quantities on the neighboring element.
At the boundary, we adopt the mirror principle following~\citet{Hesthaven2007} wherever applicable to define an outer state $(\cdot)_+$ that enforces boundary conditions weakly.

\subsection{Pressure Poisson Step}
The symmetric interior penalty Galerkin (SIPG) discretization~\cite{Arnold1982, Arnold2000, Arnold2002} is chosen for the PPE~\eqref{eq:ppe_time_discretized_leray} to avoid introducing auxiliary variables. The primal formulation reads: Find $\hat{p}^{n+1} \in \mathcal{V}_h^p$ such that
\begin{equation}
    \label{eq:fully_discrete_PPE}
    l_p\left(\hat{p}^{n+1}, q\right) 
    = r_p \left(\ve{u}^{n+1-i}_{i=1,\dots,J}, \ve{f}^{n+1}, \ve{g}^{n+1-i}_{i=1,\dots,J}, q\right)
    ,
\end{equation}
for all $q\in\mathcal{V}_h^p$ and all elements $\Omega_e \in \Omega_h$. Here, the left-hand side $l_p$ is defined as
\begin{equation}
	\label{eq:pressPoissonWeak_V1}
    \begin{aligned}
    	l_p(\hat{p}^{n+1},q) 
        \coloneqq
        &\left(\nabla \hat{p}^{n+1}, \nabla q\right)_{\Omega_e} \\
        &- \left(\nabla q, \frac{1}{2}\jump{\hat{p}^{n+1}}\ve{n}\right)_{\partial\Omega_e\backslash\Gamma_h} - \left(\nabla q, \hat{p}_-^{n+1}\ve{n}\right)_{\partial\Omega_e\cap\Gamma_N} \\
        &-\left(q, \central{\nabla \hat{p}^{n+1}}\cdot\ve{n}\right)_{\partial\Omega_e\backslash\Gamma_h}  - \left(q, \nabla \hat{p}_-^{n+1}\cdot\ve{n}\right)_{\partial\Omega_e\cap\Gamma_N} \\
        &+\left(q, \tau\jump{\hat{p}^{n+1}}\right)_{\partial\Omega_e\backslash\Gamma_h}  +\left(q, 2\tau \hat{p}_-^{n+1}\right)_{\partial\Omega_e\cap\Gamma_N},
    \end{aligned}
\end{equation}
with the stabilization parameter $\tau$ chosen according to~\citet{Hillewaert2013}.

The right-hand side $r_p$ consists of the contributions from the body force $\ve{f}$, the terms stemming from the convective operator, the time derivative, the viscous term, the Leray projection and the inhomogeneous part of the SIPG discretization:
\begin{gather}
    r_p 
    \left(\ve{u}^{n+1-i}_{i=1,\dots,J}, \ve{f}^{n+1}, \ve{g}^{n+1-i}_{i=1,\dots,J}, q\right)
    %&
    \coloneqq 
    r_\text{f}(\ve{f}^{n+1},q) 
    + \sum_{i=1}^{J_c} \Bar{\beta_i} 
     r_\text{conv}(\ve{u}^{n+1-i},\ve{g}^{n+1-i},q) 
     \nonumber \\
    %&
    \phantom{\coloneqq}
    + r_\text{SIPG}({g}_p^{n+1},\ve{u}^{n+1-i},\ve{g}^{n+1},q)
    - \sum_{i=1}^{J} \frac{\alpha_i}{\Delta t} r_\text{Leray}(\ve{u}^{n+1-i},\ve{g}^{n+1-i},q).
\end{gather}

By choosing the central flux and integrating the forcing term, the convective operator and the Leray projection term by parts once, the consistent boundary condition given in Equation~\eqref{eq:ppe_dirichlet_time_discretized} reduces to contributions from only the acceleration term and the viscous term. All remaining terms cancel on the boundary, as shown in the following. The forcing term is written as 
\begin{equation}
	r_\text{f}(\ve{f},q) 
    \coloneqq 
    (\ve{f}, \nabla q)_{\Omega_e} 
    - (\ve{f}\cdot\ve{n}, q)_{\partial\Omega_e\backslash\Gamma_h} 
    - \left(\ve{f}\cdot \ve{n} , q\right)_{\partial\Omega_e\cap\Gamma_N}.
\end{equation}
On the Dirichlet boundary the forcing term $- \left(\ve{f}\cdot \ve{n} , q\right)_{\partial\Omega_e\cap\Gamma_D}$ cancels with the contribution of the forcing term $\left(\ve{f}\cdot \ve{n} , q\right)_{\partial\Omega_e\cap\Gamma_D}$ in the consistent boundary condition. \newline
The convective term is formulated in convective formulation and is integrated by parts once. On the Dirichlet boundary, the contribution from the average operator is $\convectiveterm{\ve{u}}{\ve{g}}$ and the convective contribution from Equation~\eqref{eq:ppe_dirichlet_time_discretized} is $-\convectiveterm{\ve{u}}{\ve{g}}$, canceling each other, leading to
\begin{equation}\label{eq:rhs_ppe_convective}
\begin{aligned}
		r_\text{conv}(\ve{u},\ve{g},q) \coloneqq 
        & -\left(\convectiveterm{ \ve{u}}{\ve{u}}, \nabla q\right)_{\Omega_e} 
        %\\
		%& 
        + \left(q, \central{\convectiveterm{ \ve{u}} {\ve{u}}} \cdot\ve{n}\right)_{\partial\Omega_e\backslash\Gamma_h} 
       % \\
	%	& 
    %    + \left(q, \convectiveterm{ \ve{u}}{(\ve{g} - \ve{u})}\cdot\ve{n}\right)_{\partial\Omega_e\cap\Gamma_D} 
        \\
		& 
        + \left(q, \convectiveterm{ \ve{u}}{\ve{u}}\cdot\ve{n}\right)_{\partial\Omega_e\cap\Gamma_N}.
\end{aligned}
\end{equation}

The divergence operator for the Leray projection terms yields with boundary conditions
\begin{equation}
\begin{aligned}
       r_\text{Leray}(\ve{u},\ve{g},q) 
       \coloneqq
       %& 
       -\left(\nabla q, \ve{u} \right)_{\Omega_e} 
       + \left( q, \central{\ve{u}} \cdot \ve{n}\right)_{\partial\Omega_e\backslash\Gamma_h} 
       %\\
       %&
       +  \left(q, \ve{u} \cdot \ve{n} \right)_{\partial\Omega_e\cap\Gamma_N}.
\end{aligned}	
\end{equation}

The inhomogeneous part of the SIPG formulation includes the remaining parts of the consistent boundary condition (Equation~\eqref{eq:ppe_dirichlet_time_discretized}) and is given by
\begin{equation}
\begin{aligned}
		r_\text{SIPG}({g}_p^{n+1},\ve{u}^{n+1-i},\ve{g}^{n+1},q) \coloneqq &  
        \left(q, 2 \tau g_p^{n+1} \right)_{\partial\Omega_e\cap\Gamma_N} 
        - \left(\nabla q,  g_p^{n+1} \ve{n}\right)_{\partial\Omega_e\cap\Gamma_N} 
        \\
        &- \left( q, \frac{\gamma_0}{\Delta t} \ve{g}^{n+1} \cdot  \ve{n} \right)_{\partial\Omega_e\cap\Gamma_D}
        \\
        & - \sum_{i=1}^{J_p} \hat{\beta_i}\left(q, \nu\nabla\times\left(\nabla\times\ve{u}^{n+1-i}\right)\cdot\ve{n}\right)_{\partial\Omega_e\cap\Gamma_D},
\end{aligned}
\end{equation}
with $g_p^{n+1}$ following from Equation~\eqref{eq:bc_traction_pressure_time_discrete_leray}. Note that the second part of the time derivative involving $\sum_i^{J} \nicefrac{\alpha_i}{\Delta t}\ve{g}^{n+1-i}$ cancels with the contribution of the Leray projection on the Dirichlet boundary.
In the DG setting, the curl-curl term expressing the viscous contribution on the boundary cannot be reformulated into an expression with first derivatives only as provided in~\cite{Pacheco2021_a,Creff2025,Liu2009} for continuous elements.
Therefore, employing equal-order interpolation may lead to a loss of accuracy for highly viscous problems when adopting low polynomial degrees. Instead, our approach following \citet{Krank2016} is to first compute the vorticity $\ve{\omega} \coloneqq \nabla \times \ve{u}$ by standard element-wise $L^2$ projection to obtain the vorticity, of which we then compute the curl once again to obtain the required $\nabla \times \ve{\omega} \equiv \nabla \times \nabla \times \ve{u}$.

\subsection{Momentum Step}
The momentum step consists of finding $\ve{u}^{n+1} \in \mathcal{V}_h^u$ such that
\begin{equation}
    \label{eq:fully_discrete_momemtum_step}
    l_{u}\left(\ve{u}^{n+1}, \ve{u}^\star, \ve{g}^{n+1}, \ve{v}\right)
    =
    r_u\left(\ve{u}^{n+1-i}_{i=1,\dots,J}, \ve{u}^\star, \hat{p}^{n+1}, \ve{g}^{n+1}, g_p^{n+1}, h_u^{n+1}, \ve{v}\right)
\end{equation}
holds for all $\ve{v}\in\mathcal{V}_h^u$ and all elements $\Omega_e \in \Omega_h$.

The left-hand side $l_{u}(\ve{u}^{n+1},\ve{u}^\star,\ve{v})$ consists of the time derivative term $m(\ve{u},\ve{v})$, the convective term $c(\ve{u}, \ve{u}^\star, \ve{v})$, and the viscous part $v(\ve{u},\ve{v})$. For stabilization, a divergence penalty term $a_\text{div}(\ve{u}, \ve{u}^\star,\ve{v})$ and a penalty term $a_\text{cont}(\ve{u}, \ve{u}^\star,\ve{g},\ve{v})$ enforcing normal continuity~\cite{Joshi2016,Fehn2018b} are considered, leading to
\begin{equation}
    \begin{aligned}
        l_{u}\left(\ve{u}^{n+1}, \ve{u}^\star, \ve{g}^{n+1}, \ve{v}\right)
        &\coloneqq
        m(\ve{u}^{n+1}, \ve{v}) 
        + c(\ve{u}^{n+1}, \ve{u}^\star, \ve{v})
        + v(\ve{u}^{n+1}, \ve{v}) 
        \\
        &\phantom{\coloneqq}
        + a_\text{div}(\ve{u}^{n+1}, \ve{u}^\star, \ve{v}) 
        + a_\text{cont}(\ve{u}^{n+1}, \ve{u}^\star, \ve{g}^{n+1}, \ve{v})
        .
    \end{aligned}
\end{equation}
The time derivative term is given as
\begin{equation}
	m(\ve{u}, \ve{v}) 
    \coloneqq 
    \left(\frac{\gamma_0}{\Delta t}\ve{u}, \ve{v}\right)_{\Omega_e}.
\end{equation}
The convective term $c(\ve{u}, \ve{u}^\star, \ve{v})$ can either be formulated in the convective form $c_\text{conv}(\ve{u}, \ve{u}^\star, \ve{v})$ or the divergence form $c_\text{div}(\ve{u}, \ve{u}^\star, \ve{v})$. The upwind flux is used for the convective formulation of the convective term
\begin{equation}
	\begin{aligned}
		c_\text{conv}(\ve{u}, \ve{u}^\star, \ve{v}) 
        \coloneqq& 
        \left(\ve{u}^{\star} \cdot \nabla \ve{u}, \ve{v}\right)_{\Omega_e} 
        - 
        \left(
            \left(\central{\ve{u}^{\star}}\cdot \ve{n}\right) \ve{u}, \ve{v}
        \right)_{\partial\Omega_e\backslash\Gamma_h} 
        \\
        &+ 
        \left(
            \left(
                \central{\ve{u}^{\star}}\cdot \ve{n}
            \right)
            \central{\ve{u}}
            +
            \nicefrac{1}{2}
            \abs{\central{\ve{u}^{\star}}\cdot \ve{n}}
            \jump{\ve{u}}, \ve{v}
        \right)_{\partial\Omega_e\backslash\Gamma_h} 
        \\	
		&+ 
        \left(
            \ve{v}, 
            \left(
                \abs{\ve{u}^{\star}\cdot\ve{n}} - \ve{u}^{\star}\cdot\ve{n}
            \right)\ve{u}
        \right)_{\partial\Omega_e\cap\Gamma_D}
        .
    \end{aligned}
\end{equation}
The local Lax--Friedrichs flux is adopted within the divergence formulation of the convective term
\begin{equation}
	\begin{aligned}
		c_\text{div}(\ve{u}, \ve{u}^\star, \ve{v}) 
        \coloneqq & 
        - 
        \left(
            \ve{u} \otimes \ve{u}^{\star}, \nabla \ve{v}
        \right)_{\Omega_e} 
        \\
		&+ 
        \left(
            (\central{\ve{u} \otimes \ve{u}^{\star}} \cdot \ve{n}) 
            + 
            \nicefrac{1}{2}
            \Lambda 
            \jump{\ve{u}}, \ve{v}
        \right)_{\partial\Omega_e\backslash\Gamma_h}
        \\	
		& + 
        \left(
            \ve{v}, |\ve{u}^{\star} \cdot \ve{n}| \ve{u}
        \right)_{\partial\Omega_e\cap\Gamma_D} 
        +  
        \left(
            \ve{v}, (\ve{u} \otimes \ve{u}^{\star} ) \cdot \ve{n})
        \right)_{\partial\Omega_e\cap\Gamma_N}
        .
    \end{aligned}
\end{equation}
with $\Lambda = \zeta_\text{LF} \max (|\ve{u}^{\star}_+ \cdot \ve{n}|, |\ve{u}^{\star}_- \cdot \ve{n}|)$, where $\zeta_\text{LF}$ can be set to $0.5$ to allow for larger time steps.
\begin{remark}
    Expressing the convective operator in skew-symmetric form does not yield the same benefits (e.g., in form of energy stability) as in the continuous Galerkin setting. This difference arises from the face contributions that appear in both the convective and divergence formulations, even when choosing the same fluxes.
\end{remark}

Similar to the pressure operator, the SIPG method is utilized for the viscous part
\begin{equation}
	\begin{aligned}
        v(\ve{u}, \ve{v}) 
        \coloneqq& 
        \left(
            \nu\nabla \ve{u}, \nabla \ve{v}
        \right)_{\Omega_e} 
        - 
        \left(
            \nabla \ve{v}\cdot\ve{n}, \nicefrac{1}{2} \, \, \nu\jump{\ve{u}}
        \right)_{\partial\Omega_e\backslash\Gamma_h}
        - 
        \left(
            \ve{v}, 
            \nu
            \left(
                \central{\nabla \ve{u}\cdot\ve{n}}-\tau\jump{\ve{u}}
            \right)
        \right)_{\partial\Omega_e\backslash\Gamma_h}
        \\
        &
         -
        \left(
            \nabla\ve{v}\cdot \ve{n}, \nu\ve{u}
        \right)_{\partial\Omega_e\cap\Gamma_D} 
        - 
        \left(
            \ve{v},\nu \nabla\ve{u}
            \cdot \ve{n} - 2 \nu \tau\ve{u}
        \right)_{\partial\Omega_e\cap\Gamma_D}.
    \end{aligned}
\end{equation}
The stabilization parameter $\tau$ for tensor-product finite elements is again chosen according to~\citet{Hillewaert2013}.
The divergence penalty term is defined as
\begin{equation}
	\begin{aligned}
        a_\text{div}(\ve{u}, \ve{u}^\star, \ve{v}) 
        \coloneqq & 
        \left(
            \frac{
            \zeta_D
            h_e
            \|\overline{\ve{u}^{\star}}\| 
            }{k_u + 1}
            \,
            \nabla \cdot \ve{u}, 
            \nabla \cdot \ve{v}
        \right)_{\Omega_e}
        ,
    \end{aligned}
\end{equation}
with penalty parameter $\zeta_D = 1$, $\|\overline{\ve{u}^\star}\|$ being the element-wise volume averaged extrapolated velocity, $h_e = V_e^{1/3}$ being the characteristic element length, and with $V_e$ being the element volume~\cite{Fehn2018b}. Note that the velocity divergence penalty is related to grad-div stabilization for continuous elements~\cite{Akbas2018,Olshanski2009}.

The continuity penalty term, applied in normal direction on the inner faces and the Dirichlet boundary following~\cite{Joshi2016,Krank2016,Fehn2018b}, is defined as
\begin{equation}
	\begin{aligned}
        a_\text{cont}(\ve{u}, \ve{u}^\star, \ve{g}, \ve{v}) 
        \coloneqq &
        \left(
            \ve{v}\cdot \ve{n}, 
            \zeta_C 
            \|\overline{\ve{u}^{\star}}\| 
            (\jump{\ve{u}} \cdot \ve{n}) 
        \right)_{\partial\Omega_e\backslash\Gamma_h} 
        %\\
		%&
        + \left( 
            \ve{v}\cdot \ve{n}, 
            2 \zeta_C 
            \|\overline{\ve{u}^{\star}}\| 
            (\ve{u} - \ve{g}) \cdot \ve{n} 
        \right)_{\partial\Omega_e\cap\Gamma_D}
        ,
	\end{aligned}
\end{equation}
with penalty parameter $\zeta_C = 1$.
As the splitting method comes with consistent boundary conditions, the penalty terms can directly be incorporated in the momentum equation. It is also possible to include the terms in a postprocessing (or projection) step as is done in~\cite{Fehn2018b}.
\begin{remark}
    The divergence and continuity penalty can be seen to enforce a solution space akin to H(div) elements, where the former regulates the divergence-freeness of the velocity field and the latter the continuity in normal direction \cite{Akbas2018}.
\end{remark}

The right-hand side consists of the forcing term $f(\ve{f},\ve{v})$, the time-derivative term $m_r(\ve{u}^{n+1-i}_{i=1,\dots,J}, \ve{v})$, the pressure term $b_p(p,g_p, \ve{v})$, the convective term $c_r(\ve{u}^\star,\ve{g},\ve{v})$ and the viscous term $v_r(\ve{g}, \ve{h}_u, \ve{v})$, and reads as
\begin{equation}
    r_{u} 
    \coloneqq
    f(\ve{f}, \ve{v}) 
    + m_r(\ve{u}^{n+1-i}_{i=1,\dots,J}, \ve{v}) 
    + b_p(\hat{p}^{n+1}, g_p^{n+1}, \ve{v}) 
    + c_r(\ve{u}^\star, \ve{g}^{n+1} , \ve{v}) 
    + v_r(\ve{g}^{n+1}, \ve{h}^{n+1}_u, \ve{v})
\end{equation}
with
\begin{equation}
	\begin{aligned}
	f(\ve{f}, \ve{v}) &\coloneqq (\ve{f}^{n+1}, \ve{v})_{\Omega_e} 
    \end{aligned}
\end{equation}
and the previous time velocities of the BDF integrator,
\begin{equation}
	\begin{aligned}
        m_r(\ve{u}^{n+1-i}_{i=1,\dots,J}, \ve{v}) 
        \coloneqq &
        \left(
            \sum_{i=1}^{J}\frac{\alpha_i}{\Delta t}\ve{u}^{n+1-i}, \ve{v}
        \right)_{\Omega_e}.
    \end{aligned}
\end{equation}

The formulation for the pressure gradient term is obtained by integrating by parts twice,
\begin{equation} \label{eq:pressure_gradient_rhs}
	\begin{aligned}
    	b_p(p, g_p, \ve{v}) 
        \coloneqq& 
        - 
        \left(
            \nabla p, \ve{v}
        \right)_{\Omega_e} 
        +
        \left(
            \nicefrac{1}{2}\jump{p}\ve{n}, \ve{v}
        \right)_{\partial\Omega_e\backslash\Gamma_h}  
        + 
        \left(
            (p - g_p) \ve{n},\ve{v}
        \right)_{\partial\Omega_e\cap\Gamma_N}.
    \end{aligned}
\end{equation}
Despite the pressure gradient term only appearing on the right-hand side of the momentum balance equation, numerical results from \citet{Fehn2017} indicate that integration by parts is indeed necessary to achieve the expected convergence rates.

The inhomogeneous contributions of the convective term in convective formulation are expressed as
\begin{equation}
	\begin{aligned}
        c_{r,\text{conv}}(\ve{u}^\star, \ve{g}, \ve{v}) 
        \coloneqq 
        %& 
        -
        \left(
            \ve{v}, \left(\ve{u}^\star\cdot\ve{n}-\abs{\ve{u}^\star\cdot\ve{n}}\right)\ve{g}
        \right)_{\partial\Omega_e\cap\Gamma_D}.
    \end{aligned}
\end{equation}
For the divergence formulation the right-hand side is
\begin{equation}
	\begin{aligned}
        c_{r,\text{div}}(\ve{u}^\star, \ve{g}, \ve{v}) 
        \coloneqq 
        %& 
        \left(
            \ve{v}, -(\ve{g} \otimes \ve{u}^\star) \cdot \ve{n} + |\ve{u}^\star \cdot \ve{n}| \ve{g}
        \right)_{\partial\Omega_e\cap\Gamma_D}.
    \end{aligned}
\end{equation}
The contributions from the viscous part are due to the SIPG discretization
\begin{equation} \label{eq:viscous_sipg_rhs}
	\begin{aligned}
	v_r(\ve{g}, \ve{h}_u, \ve{v}) 
    \coloneqq & 
    - 
    \left(
        \nabla \ve{v} \ve{n}, \nu \ve{g}
    \right)_{\partial\Omega_e\cap\Gamma_D} 
    + 
    \left(
        \ve{v}, 2\nu\tau\ve{g}
    \right)_{\partial\Omega_e\cap\Gamma_D}
    %\\
	%&
    +
    \left(
        \ve{h_u},\ve{v}
    \right)_{\partial\Omega_e\cap\Gamma_N},
	\end{aligned}
\end{equation}

with $\ve{h}_u = \nu (\nabla \ve{u} )\ve{n}$. Note that combining Equations~\eqref{eq:pressure_gradient_rhs} and \eqref{eq:viscous_sipg_rhs} recovers the complete traction $\ve{h} = \nu (\nabla \ve{u} )\ve{n} - g_p \ve{n}$, such that only $\ve{h}$ must be specified on the Neumann boundary.

\subsection{Timestepping Algorithm}
This section summarizes the time stepping procedure of the modified pressure scheme in Alg.~\ref{alg:consistent_splitting_scheme}. For simplicity of presentation, the startup process is omitted here, which either consists of interpolating a given exact solution at the required time instants, using lower-order BDF schemes with increasing order, or other time integrators of appropriate order until the necessary data $\ve{u}^{n+1-i}$ is gathered.

\begin{algorithm}
	\begin{algorithmic}[1]
		\caption{Splitting solver}
		\label{alg:consistent_splitting_scheme}
		\Function{AdvanceSplittingSolver}{$\ve{u}^{n+1-i}_{i=1,\dots,J}$, $\Delta t$}
        \State \textbf{parameter update:} select/update time integration and extrapolation constants
        \State \textbf{pressure step:} compute the pressure $\hat{p}^{n+1}$
		via~\eqref{eq:fully_discrete_PPE}, i.e.,
        find $\hat{p}^{n+1} \in \mathcal{V}_h^p$ such that
        \begin{align*}
            l_p(\hat{p}^{n+1}, q) 
            = r_p (\ve{u}^{n+1-i}_{i=1,\dots,J}, \ve{f}^{n+1}, \ve{g}^{n+1-i}_{i=1,\dots,J}, q)
            \quad \forall q\in\mathcal{V}_h^p, \forall \Omega_e \in \Omega_h
            .
            \\[-2mm]
        \end{align*}
		\State compute velocity extrapolation $\ve{u}^\star$ for linearization
		\State \textbf{momentum step:} compute the velocity
		$\ve{u}^{n+1}$ via~\eqref{eq:fully_discrete_momemtum_step}, that is, find
        $\ve{u}^{n+1} \in \mathcal{V}_h^u$ such that

        \begin{align*}
            &l_{u}(\ve{u}^{n+1}, \ve{u}^\star, \ve{g}^{n+1}, \ve{v})
            \\
            &=
            r_u(\ve{u}^{n+1-i}_{i=1,\dots,J}, \ve{u}^\star, \hat{p}^{n+1}, \ve{g}^{n+1}, g_p^{n+1}, h_u^{n+1}, \ve{v})
            \quad \forall \ve{v}\in\mathcal{V}_h^u, \forall \Omega_e \in \Omega_h
            .
            \\[-2mm]
        \end{align*}

		\State \Return $\ve{u}_h^{n+1}$, $p_h^{n+1}$
		\EndFunction
	\end{algorithmic}
\end{algorithm}

\section{Numerical Results}
\label{sec:numerical_results}

This section examines the numerical properties of the devised method. First, convergence in space and time for an example with known analytical solution is demonstrated in Sec.~\ref{sec:numerical_results_manufatured}, thereafter, the schemes' ability to accurately predict fluid mechanical quantities in practically relevant scenarios is evaluated for classical benchmarks being the two-dimensional flow around a cylinder~\cite{Schaefer1996} and the Taylor--Green vortex~\cite{TaylorGreen1937} in its three-dimensional setup according to~\citet{Wang2013}. Matrix-based and matrix-free implementations are publicly available via the \texttt{BDFConsistentSplitting} scheme in~\texttt{ExaDG}\footnote{\url{https://github.com/exadg/exadg}, retrieved on November 30, 2025.}~\cite{exadg-2020},
which is based on~\texttt{deal.II}~\cite{dealii}. The PPE is solved using a conjugate-gradient solver preconditioned by multigrid algorithms, while the momentum equation is solved with a GMRES solver using an inverse mass preconditioner, as the resulting linear system is not symmetric. If not specified otherwise, the relative solver tolerance is set to $10^{-6}$ and the absolute solver tolerance to $10^{-12}$, measured in the unpreconditioned residual. Results are presented for the modified pressure scheme  unless explicitly noted otherwise.

\subsection{Manufactured Solution}
\label{sec:numerical_results_manufatured}
First, we focus on the two dimensional Taylor--Green vortex test case~\cite{TaylorGreen1937}, see also~\citet{Hesthaven2007}. The solution of the incompressible Navier--Stokes problem is given by
\begin{equation}
    \begin{aligned}
        \ve{u}(\ve{x},t) &= \binom{-\sin(2\pi x_2)}{\phantom{+}\sin(2\pi x_1)} \exp(-4\nu \pi^2t ) \\
        p(\ve{x},t) &= -\cos(2\pi x_1) \cos(2\pi x_2) \exp(-8\nu \pi^2t)
    \end{aligned}
\end{equation}
with the right-hand side being zero, $\ve{f} = \ve{0}$. The domain $\Omega \coloneqq [-0.5,0.5]^2$ , a square of length 1, is considered and the simulation is performed over the time interval $0\leq t \leq T = 1$ with the viscosity $\nu =  0.025$. The Dirichlet and Neumann boundary conditions and initial conditions are derived from the exact solution. The relative $L^2$ errors of the velocity and the pressure are computed as 
\begin{equation}
    \epsilon_p = \left( \frac{||p-p_h||_{L^2(\Omega)}}{||p||_{L^2(\Omega)}} \right) \bigg\rvert_{t=T} \text{ and } \epsilon_u = \left( \frac{||\ve{u}-\ve{u}_h||_{L^2(\Omega)}}{||\ve{u}||_{L^2(\Omega)}} \right) \bigg\rvert_{t=T}
\end{equation}

\subsubsection{Convergence in Time}
\begin{figure}[!h]
    \centering
    \includegraphics[width=0.99\linewidth]{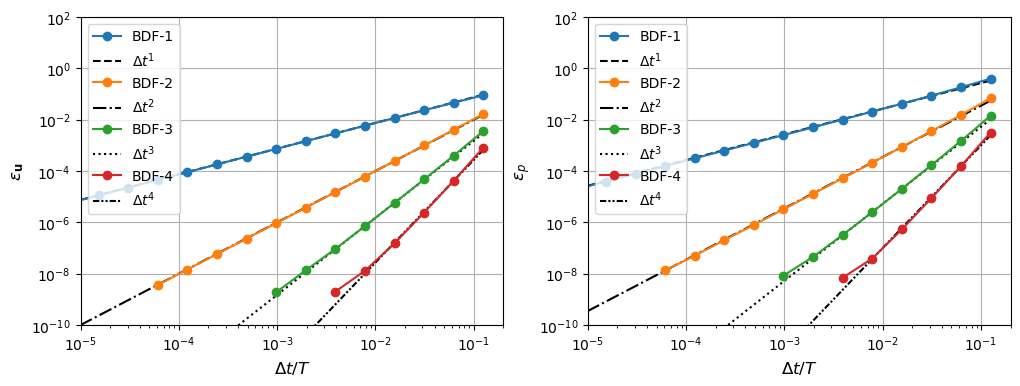}
    \caption{Temporal convergence of velocity (left) and pressure (right) in the 2D Taylor--Green vortex benchmark as obtained via the consistent splitting scheme for BDF-1 to BDF-4.}
    \label{fig:temporal_convergence}
\end{figure}

To verify the expected convergence rates in time, the grid is uniformly refined four times, yielding a uniform Cartesian grid with $h_e=1/32$. Quadrilateral elements with polynomial degree $k_{u}=5$ for the velocity and $k_p=4$ for the pressure are used. Therefore, the error in the spatial discretization is small compared to the error in time. For BDF-1 and BDF-2 $J_c$ and $J_p$ are chosen equal to $J$, while for the other cases, $J_p = J_c = J-1$ is used. 
Fig.~\ref{fig:temporal_convergence} shows optimal convergence for BDF-1 to BDF-4, that is, $\mathcal{O}(\Delta t^{J})$, for both the velocity and the pressure \textit{even though} $J_c = J-1$ for $J>2$. In the present scenario, one reaches the level of the spatial error, without observable evidence of a reduction in the convergence rate for $J>2$. However, a reduction to an order of $J_c = J-1$ cannot be excluded in general. A thorough analysis of this potential behavior falls outside the scope of the current study.

\subsubsection{Convergence in Space}
\begin{figure}[h!]
    \centering
   \begin{subfigure}{0.5\textwidth}
  \centering
  \includegraphics[width=0.99\linewidth]{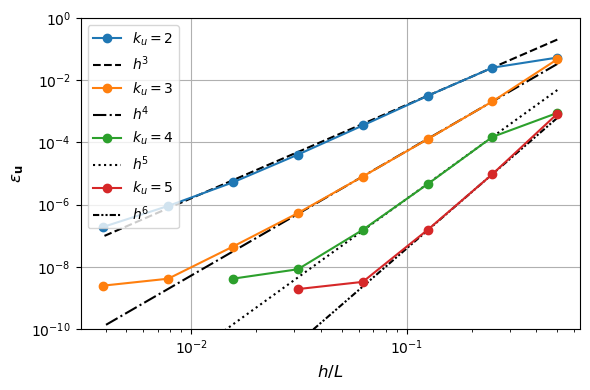}
  \caption{Convergence of the velocity}
\end{subfigure}%
\begin{subfigure}{.5\textwidth}
  \centering
  \includegraphics[width=0.99\linewidth]{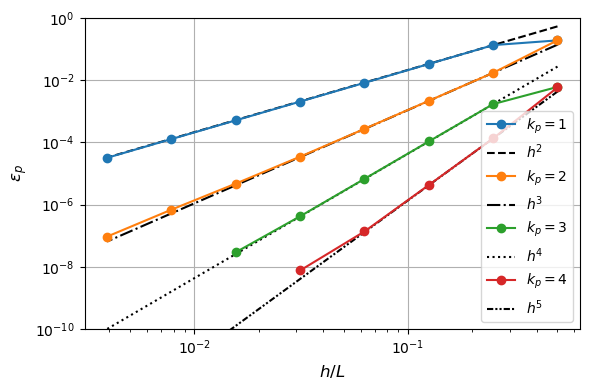}
  \caption{Convergence of the pressure}
\end{subfigure}
    \caption{Spatial convergence of the consistent splitting scheme with different polynomial degrees.}
    \label{fig:spatial_convergence}
\end{figure}

To compute the spatial error, the time step size is set to $6.1035\times10^{-5}$ applying the BDF-2 scheme, such that the temporal error is lower than the spatial one. The relative $L^2$ error under mesh refinement is depicted in Fig.~\ref{fig:spatial_convergence}. Optimal convergence rates of $h^{k+1}$ can be observed for all polynomial degrees $k_{u}=2,\ldots,5$ and $k_p=1,\ldots,4$ for the velocity and the pressure, respectively. For $k_u = 2$ the absolute tolerance of the linear solvers is set to $10^{-8}$ to achieve proper convergence rates for the two finest grids.

\begin{remark}
    Adding the divergence and continuity penalty terms does not alter the order of convergence, nor affects the magnitude of the error in this specific example. However, disregarding the Leray projection---which in the present form of the splitting scheme merely leads to a modified pressure variable---decreases the convergence order of the method: instead of the optimal convergence rate of $h^{k_p+1}$, only $h^{k_p}$ is observed for the pressure, while the velocity shows rates of $h^{k_u-1}$ without the Leray projection instead of the optimal $h^{k_u+1}$. Similarly, the convergence order in time without Leray projection is impacted. This leads to the conclusion that the addition of the Leray projection is vital for the performance of the presented scheme. Implementation effort and operator complexity are marginally affected, since no additional linear system needs to be solved during the time step.
\end{remark}

\subsubsection{Comparison of the Consistent Splitting Variants}
\begin{figure}[!h]
    \centering
    \includegraphics[width=0.99\linewidth]{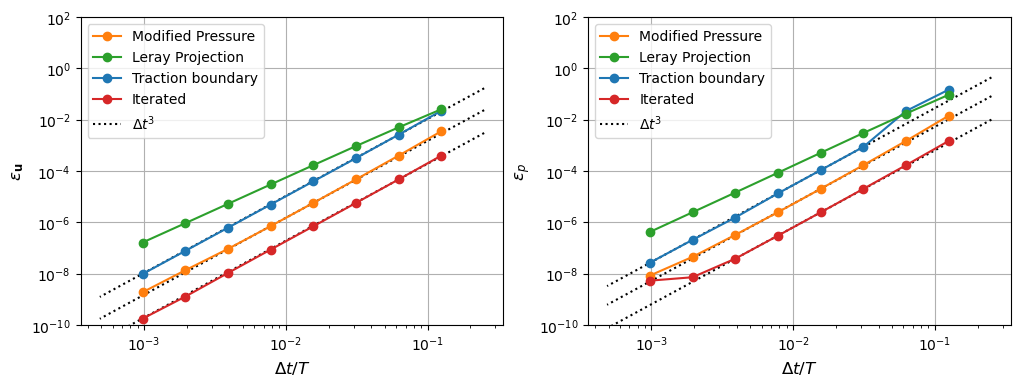}
    \caption{Comparison of the different variants of the consistent splitting scheme as summarized in Table~\ref{tab:overview_schemes} using BDF-3. The Leray corrected scheme shows suboptimal convergence rates due to the reduced extrapolation order, required for stability. For the schemes with modified pressure variable, this reduction is not observed.}
    \label{fig:temporal_convergence_bdf3}
\end{figure}

Tables~\ref{tab:consistent_splitting_variants} and \ref{tab:consistent_splitting_variants_pressure} list the relative $L^2$-errors and convergence rates obtained with the variants of the scheme presented in Chapter~\ref{chap:variants}. All variants achieve the optimal convergence rates, apart from the Leray Projection scheme. In this variant, the momentum equation is solved first with the extrapolated pressure and the Leray correction is applied in the acceleration term (cf. Table~\ref{tab:overview_schemes}). Using $J=3$ or $J=4$ for this scheme leads to some temporal instabilities, requiring $J_c = J-1$, predominantly for the latter case. This slightly lowers the achieved temporal convergence rate. A visualization of this data for BDF-3 is provided in Figure~\ref{fig:temporal_convergence_bdf3}, showing that the iterated version results in the lowest error, as it resolves the nonlinearities up to the chosen tolerance. We also provide results adopting the full traction boundary condition~\eqref{eq:bc_traction_pressure_time_discrete_leray}, denoted by ``Traction''. This requires an extrapolation of the velocity in the Neumann boundary term of the PPE, increasing the error compared to the modified pressure scheme, for which the pressure data $g_p$ is assumed known. Additionally, $J_h = J$ was decreased to $J_h = 2$ when considering the full traction vector $\ve{h}$ for stability reasons for BDF-3 and higher. Based on these presented results, the modified pressure version seems to be the best choice out of the presented variants, as it combines high computational efficiency with accurate results.

\begin{table}[ht]
\centering
\caption{Relative velocity $L^2$-errors and observed convergence rates in parentheses for BDF orders 1 to 4 and various consistent splitting scheme variants as summarized in Table~\ref{tab:overview_schemes}. For the corresponding pressure errors and convergence rates see Table~\ref{tab:consistent_splitting_variants_pressure}.}
\label{tab:consistent_splitting_variants}
\resizebox{\linewidth}{!}{
\begin{tabular}{|c|c|c|c|c|c|}
\hline
Scheme & $\Delta t/T \times10^{-2}$ &  BDF-1 & BDF-2 & BDF-3 & BDF-4 \\
\hline

Leray
& $3.13$ & 
3.46e-2 (1.15) & 
1.48e-3 (2.33) & 
9.15e-4 (2.43) & 
5.24e-5 (3.42) \\
Projection 
& $1.65$ & 
1.58e-2 (1.13) & 
3.01e-4 (2.29) & 
1.64e-4 (2.48) & 
4.88e-6 (3.42) \\
& $0.78$ & 
7.30e-3 (1.11) & 
6.37e-5 (2.24) & 
2.90e-5 (2.50) & 
4.31e-7 (3.50) \\
\hline
 
Modified 
& $3.12$ & 
2.35e-2 (0.98) & 
1.00e-3 (2.03) & 
4.77e-5 (3.08) & 
2.43e-6 (4.11) \\
Pressure
& $1.56$ & 
1.18e-2 (0.99) & 
2.47e-4 (2.02) & 
5.79e-6 (3.04) & 
1.54e-7 (3.98) \\
& $0.78$ & 
5.90e-3 (1.00) & 
6.11e-5 (2.01) & 
7.18e-7 (3.01) & 
1.31e-8 (3.55) \\
\hline

\multirow{3}{*}{Iterated} 
& $3.12$ & 
1.18e-2 (0.99) & 
2.47e-4 (2.02) & 
5.78e-6 (3.04) & 
1.46e-7 (4.05) \\
& $1.56$ & 
5.90e-3 (1.00) & 
6.11e-5 (2.01) & 
7.10e-7 (3.03) & 
9.40e-9 (3.95) \\
& $0.78$ & 
2.94e-3 (1.00) & 
1.52e-5 (2.01) & 
8.76e-8 (3.02) & 
9.42e-10 (3.32) \\
\hline

\multirow{3}{*}{Traction} 
& $3.12$ & 
2.16e-2 (1.20) &
6.94e-4 (2.50) &
3.22e-4 (3.03) &
1.84e-5 (4.27) \\
& $1.56$ & 
9.96e-3 (1.12) & 
1.41e-4 (2.30) & 
4.00e-5 (3.01) & 
1.11e-6 (4.05) \\
& $0.78$ & 
4.75e-3 (1.07) & 
3.15e-5 (2.17) & 
4.99e-6 (3.00) & 
6.81e-8 (4.03) \\
\hline

\end{tabular}}
\end{table}

\subsection{Flow Around a Cylinder}
\label{sec:numerical_results_cylinder}
\begin{figure}
    \centering
    \begin{overpic}[width=0.8\linewidth]{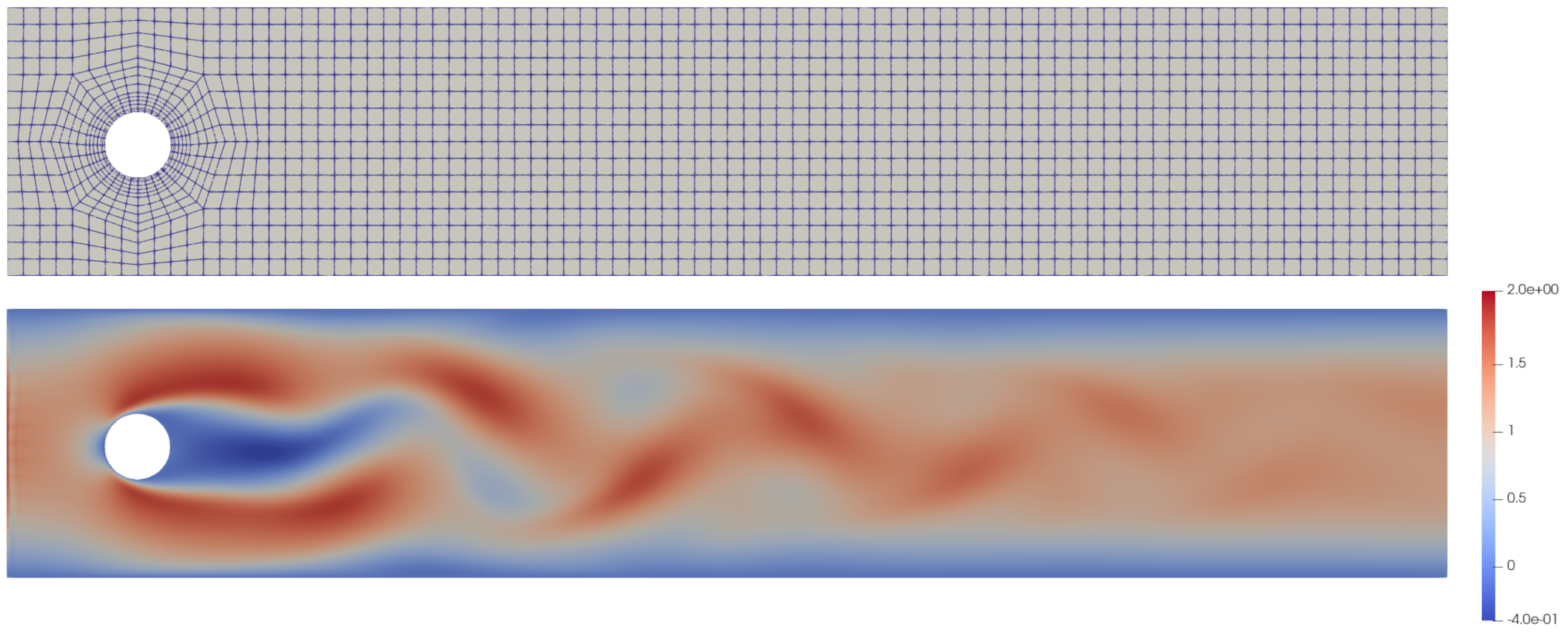}
        \put(98, 17.5){\footnotesize \rotatebox{270}{velocity $u_x$}}
    \end{overpic}
    \caption{Flow around a cylinder benchmark: computational grid after 2 uniform refinements (top) and horizontal component $u_x$ of the velocity vector with developed vortex shedding (bottom).}
    \label{fig:flow_around_cylinder}
\end{figure}

The flow around a cylinder benchmark, specifically the 2D-3 case proposed by \citet{Schaefer1996}, is investigated with the aim to demonstrate the effect of employing higher-order polynomial degrees in a practically relevant setting. The geometry considered is a channel with length $L=2.2$ and height $H=0.41$. A cylinder centered at $(0.2, 0.2)$ with a diameter of $D=0.1$ is placed in the channel, see Fig.~\ref{fig:flow_around_cylinder}. The inflow boundary condition at $x_1 = 0$ is given as
\begin{equation}
    \ve{g}_u(x_1=0, x_2, t) = \binom{4U_m \frac{x_2(H-x_2)}{H^2}\sin(\pi t/T)}{0}. 
\end{equation}
The average velocity is given as $U_m=1.5$ and the viscosity as $\nu=10^{-3}$ such that the resulting maximal Reynolds number is $Re_{max}=100$. The time interval is chosen as $0\leq t \leq T = 8$. At $x_1=L$, an open traction (outflow) boundary is given with $g_p=0$, while at all other boundaries the no-slip condition, $\ve{u} = \ve{0}$, is prescribed. 

\subsubsection{Comparison to reference values}

The results obtained via the splitting scheme using different polynomial degrees are compared with values from the literature \citep{Schaefer1996,John2004,Fehn2017}. The results in~\cite{Fehn2017} are obtained with $k_u = 10$ with around $10^6$ DoFs. Tab.~\ref{tab:cylinder_flow_reference_values} lists the reported reference values for the maximum lift and drag values and the pressure difference across the cylinder at time $t=T$. Tab.~\ref{tab:cylinder_flow_values} lists the values computed with the present method under two refinements (a total of 800 elements, similar to \cite{Fehn2017}) at different polynomial degrees, using a BDF-3 scheme with a time step of $10^{-5}$. The results show excellent agreement of all quantities for all velocity degrees considered with the literature results.

\begin{table}[ht]
  \centering
  \caption{Reference values for the maximum drag and lift values, $c_{D,\max}$ and $c_{L,\max}$, and the pressure difference $\Delta p$ at $t=T$ from literature for the 2D-3 test case of the flow around a cylinder benchmark~\cite{Schaefer1996}.}
  \label{tab:cylinder_flow_reference_values}
  \resizebox{\linewidth}{!}{
  \begin{tabular}{|l|l|l|l|}
    \hline
    Reference & \multicolumn{1}{|c|}{$c_{D,\max}$} 
              & \multicolumn{1}{|c|}{$c_{L,\max}$}
              & \multicolumn{1}{|c|}{$\Delta p(t = T)$} \\
    \hline
    \citet{Schaefer1996}~(1996) & $2.95 \pm 2 \cdot 10^{-2}$ & $0.48 \pm 1 \cdot 10^{-2}$ & $-0.11 \pm 5 \cdot 10^{-3}$ \\
    \citet{John2004}~(2004)           & $2.950921575 \pm 5 \cdot 10^{-7}$ & $0.47795 \pm 1 \cdot 10^{-4}$ & $-0.1116 \pm 1 \cdot 10^{-4}$ \\
    \citet{Fehn2017}~(2017)    & $2.95091839$ & $0.47788776$ & $-0.11161590$ \\
    \hline
  \end{tabular}}
\end{table}

\begin{table}[ht]
  \centering
  \caption{Maximum drag and lift values $c_{D,\max}$ and $c_{L,\max}$ as obtained with the present scheme for test case 2D-3 using velocity degrees $k_u = 3,4,5$ showing excellent agreement with the reference values given in Tab.~\ref{tab:cylinder_flow_reference_values}.}
  \label{tab:cylinder_flow_values}
  \begin{tabular}{|l|c|c|c|}
    \hline
    Degree $k_u$ & $c_{D,\max}$ & $c_{L,\max}$ & $\Delta p(t = T)$ \\
    \hline
    3 & $2.95096212$ & $0.4778764498$ & $-0.11163197467$ \\
    4 & $2.95092039$ & $0.4778875586$ & $-0.11161954441$  \\
    5 & $2.95091848$ & $0.4778877448$ & $-0.11161596715$ \\
    %8 & $2.95091839$ & $0.4778874054$ & $-0.11161790449$ \\
    \hline
  \end{tabular}
\end{table}

\subsubsection{Higher order time integrators and temporal stability}
\label{sec:numerical_results_stability_cylinder}

To showcase the potential of higher-order time integration with the consistent splitting scheme, we demonstrate temporal stability properties for the various linearization variants of the convective term. Tab.~\ref{tab:cylinder_flow_cfl} presents results obtained considering the convective term explicit, linearly implicit or fully implicit under increasing target CFL numbers controlling the time step size. For the consistent splitting scheme, the convective formulation of the convective term is chosen, as it permits larger time-step sizes than the divergence formulation in this example. Presented is the maximum drag coefficient obtained with $k_u = 5$ in a mesh with 800 elements. This spatial discretization is unchanged when increasing the target CFL number. The explicit method remains stable only for small CFL numbers, whereas the linearly implicit formulation remains robust for CFL numbers up to at least 16. Surprisingly, the fully implicit method does not reach the same CFL limit as the linearly implicit scheme. However, for large CFL values the implicit version yields results that are closer to the reference solution, albeit at the expense of solving a nonlinear system.

\begin{table}[ht]
  \centering
  \caption{Maximum drag values $c_{D,\text{max}}$ (reference value $2.95091839$ from \cite{Fehn2017}, see Tab.~\ref{tab:cylinder_flow_reference_values}) for test case 2D-3 using the explicit, linearly implicit and implicit versions of the convective term under different CFL numbers with BDF-2. Simulations that fail to converge are marked by a dash. }
  \label{tab:cylinder_flow_cfl}
  \resizebox{\linewidth}{!}{
  \begin{tabular}{|l|c|c|c|c|c|c|}
    \hline
    CFL number & $0.4$ & $1.0$ & $2.0$ & $4.0$ & $8.0$ & $16.0$ \\
    \hline
    Explicit & $2.95091536$ & $2.95091496$ & $-$ & $-$& $-$& $-$ \\
    Linearly implicit & $2.95091557$ & $2.95091625$ & $2.95092445$ & $2.95096155$ & $2.95113375$ & $2.95232103$ \\
    Implicit & $2.95091523$ & $2.95091335$ & $2.95091206$ & $2.95091085$ & $2.95091590$ & $-$  \\
  
    \hline
  \end{tabular}}
\end{table}

For a comparison to the dual splitting scheme, Tab.~\ref{tab:cylinder_flow_bdf} presents results obtained for the consistent splitting scheme and the dual splitting scheme in their linearly implicit formulation for different BDF schemes under increasing CFL numbers. For the dual splitting scheme, the continuity and the divergence penalty terms are added in a postprocessing step (see~\cite{Fehn2022}). For better comparison, the same is done in case of the consistent-splitting scheme. Both schemes remain unconditionally stable employing BDF schemes up to and including order 2 (i.e. choosing $J=J_c=J_p=2$), whereas differences arise for BDF-3 and higher. For the consistent splitting scheme, the extrapolation order of the convective term on the right-hand side of the PPE is set to $J_c = J-1$, which permits larger CFL numbers than the dual splitting scheme, where this term is extrapolated with order $J_c = J$. Note that in this case both schemes use $J_p = J - 1$.

Moreover, the increased accuracy of BDF-3 relative to BDF-2 becomes evident once $\text{CFL}>1.0$.

\begin{table}[ht]
  \centering
  \caption{Maximum drag values $c_{D,\text{max}}$ (reference value $2.95091839$ from \cite{Fehn2017}, see Tab.~\ref{tab:cylinder_flow_reference_values}) obtained for test case 2D-3 using the linearly implicit consistent splitting (CS) and dual splitting (DS) schemes with different BDF time integration schemes. Simulations that fail to converge are marked by a dash. }
  \label{tab:cylinder_flow_bdf}
  \resizebox{\linewidth}{!}{
  \begin{tabular}{|l|l|c|c|c|c|c|c|}
    \hline
    \multicolumn{2}{|c|}{{CFL number}} & {0.4} & {1.0} & {2.0} & {4.0} & {8.0} & {16.0} \\
    \hline
    \multirow{3}{*}{{CS}} 
      & BDF-2 & 2.95091557 & 2.95091625 & 2.95092445 & 2.95096155 & 2.95113375 & 2.95232103\\
      & BDF-3 & 2.95091538 & 2.95091365 & 2.95092369 & 2.95092490  & 2.95093743 & 2.95106423 \\
      & BDF-4 & 2.95092587 & 2.95092424 & -- & -- & -- & -- \\
    \hline
    \multirow{3}{*}{{DS}} 
      & BDF-2 & 2.95092625 & 2.95092711 & 2.95093520 & 2.95097113 & 2.95113559 & 2.95222858 \\
      & BDF-3 & 2.95092583 & 2.95092423 & -- & -- & -- & -- \\
      & BDF-4 & 2.95092587 & -- & -- & -- & -- & -- \\
    \hline
  \end{tabular}}
\end{table}

\subsection{3D Taylor--Green Vortex}
\label{sec:numerical_results_vortex}

To investigate the performance of the consistent splitting scheme for turbulent flows, we analyze the 3D Taylor--Green vortex problem~\cite{TaylorGreen1937} in the definition proposed by \citet{Wang2013}, representing transition from a smooth initial condition into turbulence. We consider the spatial domain being  a cube with edge length $2\pi$ centered at the origin, i.e., $\Omega_h = [-\pi, \pi]^3$. The initial condition is given by
\begin{equation}
    \begin{aligned}
        u_1(\ve{x},t=0) &= \sin(x_1) \cos(x_2) \cos(x_3) \\
        u_2(\ve{x},t=0) &= -\cos(x_1) \sin(x_2) \cos(x_3) \\
        u_3(\ve{x},t=0) &= 0 \\
        p(\ve{x},t=0) &= \frac{1}{16}\left( \cos(2x_1)+\cos(2x_2)\right) \left( \cos(2x_3) + 2\right)
    \end{aligned}
\end{equation}
with $0 \leq t \leq T = 20$ with time step selection aiming for $\mathrm{CFL}=0.4$. 
The viscosity is chosen such that the Reynolds number is $\mathrm{Re}=\nicefrac{1}{\nu} = 1600$. 
Periodic boundary conditions are applied to the boundary.

The quantities of interest are the total kinetic energy $E$ and the molecular energy dissipation rate~$\epsilon$ defined as (see, e.g., \cite{Gassner2013, Fehn2018b}) 
\begin{equation}
    E 
    \coloneqq 
    \frac{
        \int_\Omega \nicefrac{1}{2}\,\, \ve{u} \cdot \ve{u} \,\mathrm{d}\ve{x}
    }{
        \int_\Omega 1 \,\mathrm{d}\ve{x}
    }
    ,
    \qquad
    \epsilon
    \coloneqq
    \frac{
        \nu \int_\Omega \nabla \ve{u} : \nabla \ve{u} \,\mathrm{d}\ve{x}
    }{
        \int_\Omega 1 \,\mathrm{d}\ve{x}
    }
\end{equation}
and the numerical dissipation of the discrete scheme obtained as $-\frac{\mathrm{d}E}{\mathrm{d}t}-\epsilon$. 

\subsubsection{Comparison to reference values}

\begin{figure}[h!]
    \centering
   \begin{subfigure}{0.328\textwidth}
  \centering
  \includegraphics[width=0.99\linewidth]{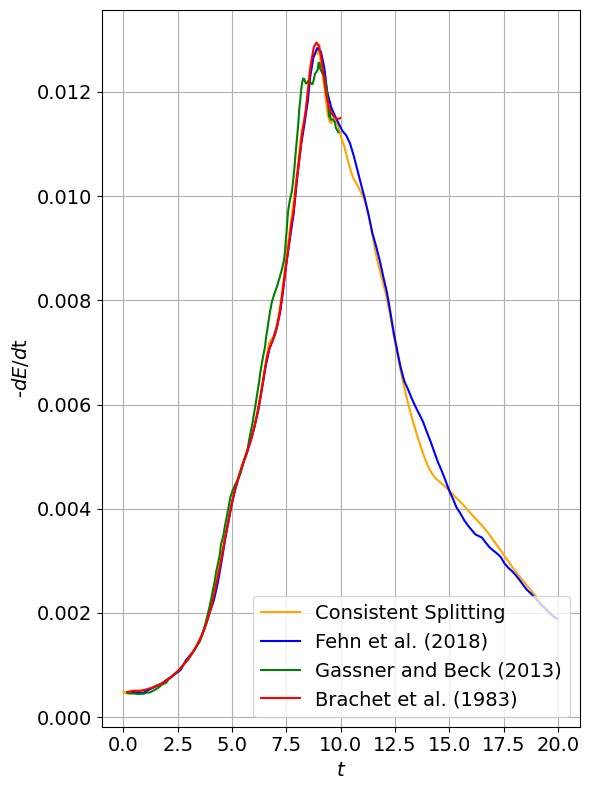}
  \caption{Kinetic energy dissipation rate}
\end{subfigure}%
\begin{subfigure}{.328\textwidth}
  \centering
  \includegraphics[width=0.99\linewidth]{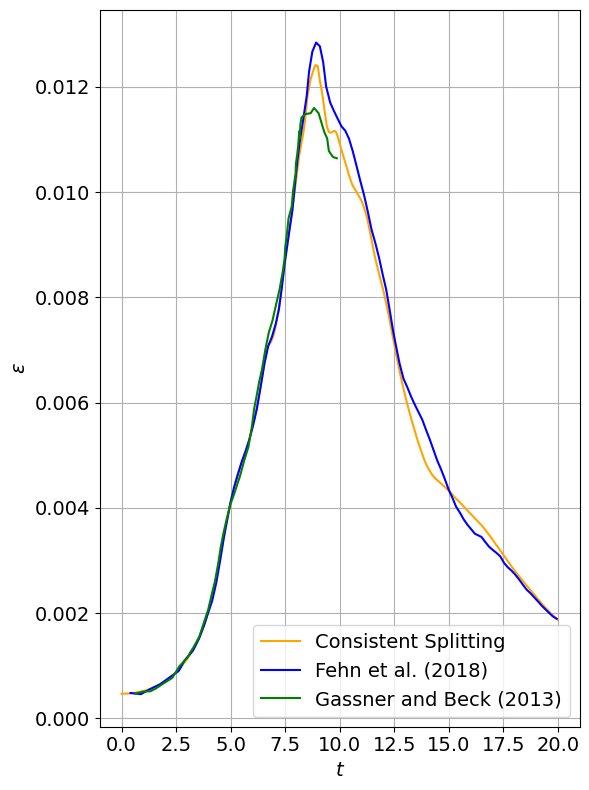}
  \caption{Molecular energy dissipation rate}
\end{subfigure}
\begin{subfigure}{.328\textwidth}
  \centering
  \includegraphics[width=0.99\linewidth]{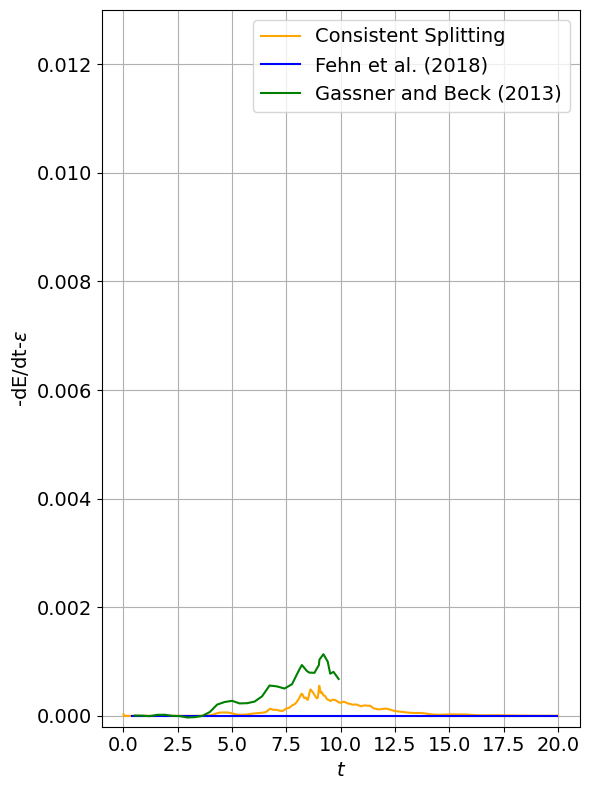}
  \caption{Numerical dissipation rate}
\end{subfigure}
    \caption{Kinetic energy, molecular dissipation and numerical dissipation of the consistent splitting method with 21 million velocity DoFs in comparison to literature results taken from \cite{Brachet1983, Gassner2013, Fehn2018b}. Values from \cite{Gassner2013} are only available for $t\in[0,10]$. Lower numerical dissipation means more accurate, limited by the fixed spatiotemporal resolution chosen.
    }
    \label{fig:taylor_green_dissipation}
\end{figure}

Comparing the consistent splitting scheme to reference values~\cite{Fehn2018b, Gassner2013, Brachet1983} in Fig.~\ref{fig:taylor_green_dissipation}, one observes that the present scheme indeed yields accurate results. Here, it should be noted that the values in \citet{Fehn2018b} are obtained via a discretization with 3.22 billion DoFs for the velocity unknowns (refinement level 8 leading to $256^3$ cells and a polynomial degree of 3) while the results using the presented splitting method were obtained using 21 million velocity DoFs (refinement level 5 yielding $32^3$ cells and a polynomial degree of 5, giving an effective resolution of $192^3$ spatial points). \citet{Brachet1983} obtained their results by direct numerical simulation and \citet{Gassner2013} utilized implicit large eddy simulation. 

\begin{figure}[h!]
    \centering
   \begin{subfigure}{0.328\textwidth}
  \centering
  \includegraphics[width=0.99\linewidth]{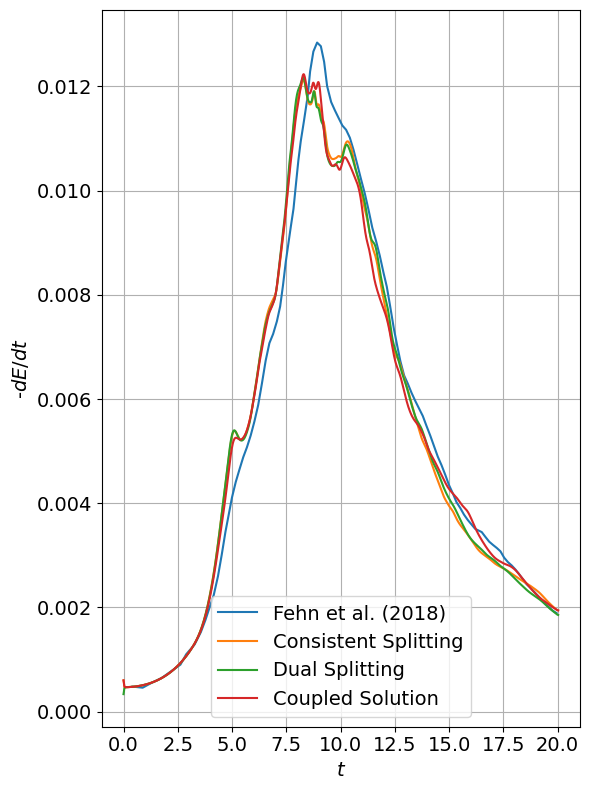}
  \caption{Kinetic energy dissipation rate}
\end{subfigure}
\begin{subfigure}{.328\textwidth}
  \centering
  \includegraphics[width=0.99\linewidth]{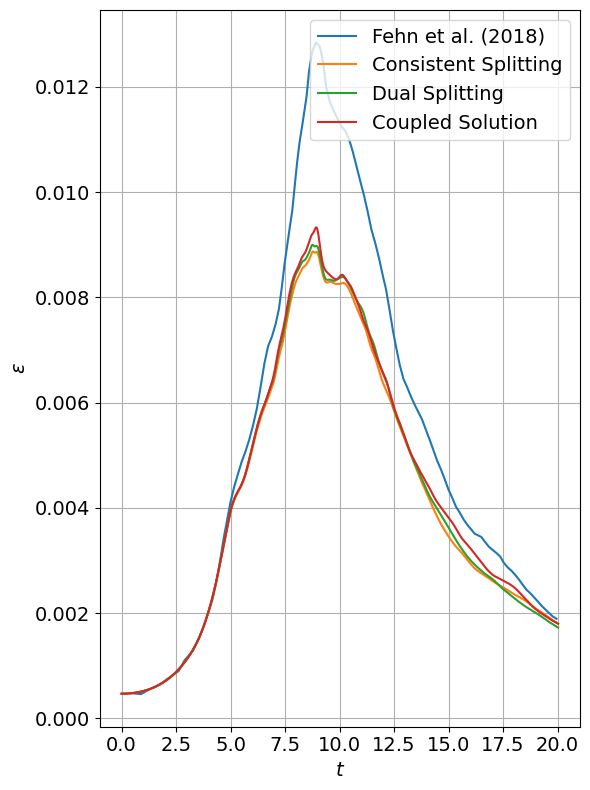}
  \caption{Molecular energy dissipation rate}
\end{subfigure}
\begin{subfigure}{.328\textwidth}
  \centering
  \includegraphics[width=0.99\linewidth]{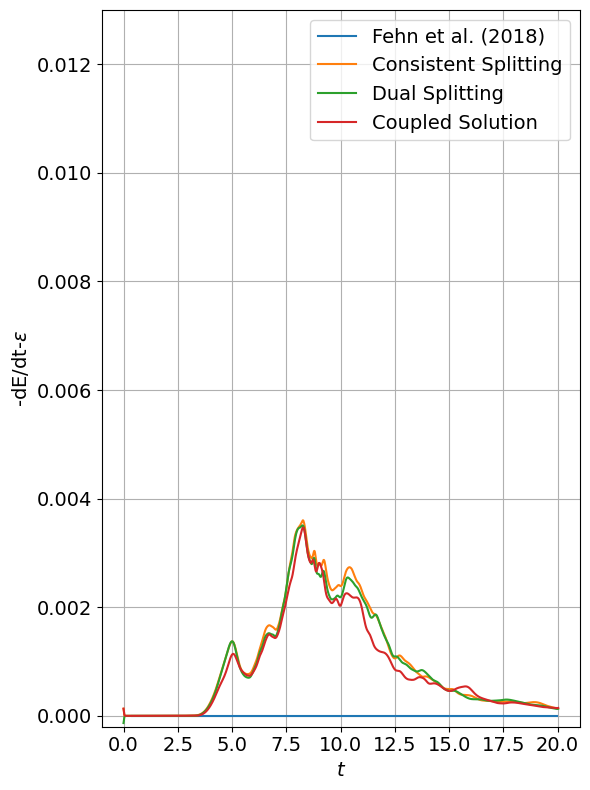}
  \caption{Numerical dissipation rate}
\end{subfigure}
    \caption{Comparison between the solutions of the coupled solution approach, the dual splitting scheme and the consistent splitting with 786432 velocity DoFs. All schemes yield similar results and compare well against the reference solution from \citet{Fehn2018b}, which uses the dual splitting scheme with 3.221 billion DoFs. Lower numerical dissipation means more accurate, limited by the fixed spatiotemporal resolution chosen.}
    \label{fig:taylor_green_dissipation_coupled_dual}
\end{figure}

For a direct comparison of splitting, projection and coupled time stepping algorithms, Fig.~\ref{fig:taylor_green_dissipation_coupled_dual} shows a comparison with identical spatial resolution with velocity degree $k_u=7$ yielding a total of 786,432 velocity DoFs and using different time stepping algorithms:
\begin{enumerate}[label=\roman*)]
    \item \textbf{dual splitting scheme}: fully explicit, CFL=0.4, convective term in convective formulation~\cite{Fehn2017},
    \item \textbf{coupled solution approach}: fully implicit, CFL=0.4, convective term in convective formulation (see, e.g., \cite{Fehn2017, Schussnig2025b} among others),
    \item \textbf{consistent splitting scheme}: linearly implicit, CFL=0.4, convective term in convective formulation (this work).
\end{enumerate}
Results demonstrate that all these methods yield comparable results, only showing minor differences between the schemes using \textit{the same spatial resolution}. For reference, the reference results from~\citet{Fehn2018b} are included once again in Fig.~\ref{fig:taylor_green_dissipation_coupled_dual}. All three schemes, i.e., the dual splitting scheme, the coupled solution approach as well as the splitting scheme are implemented in~\texttt{ExaDG}~\cite{Fehn2022, exadg-2020}. 

\subsubsection{Higher order time integrators and temporal stability}
\begin{figure}[h!]
    \centering
   \begin{subfigure}{0.33\textwidth}
  \centering
  \includegraphics[width=0.99\linewidth]{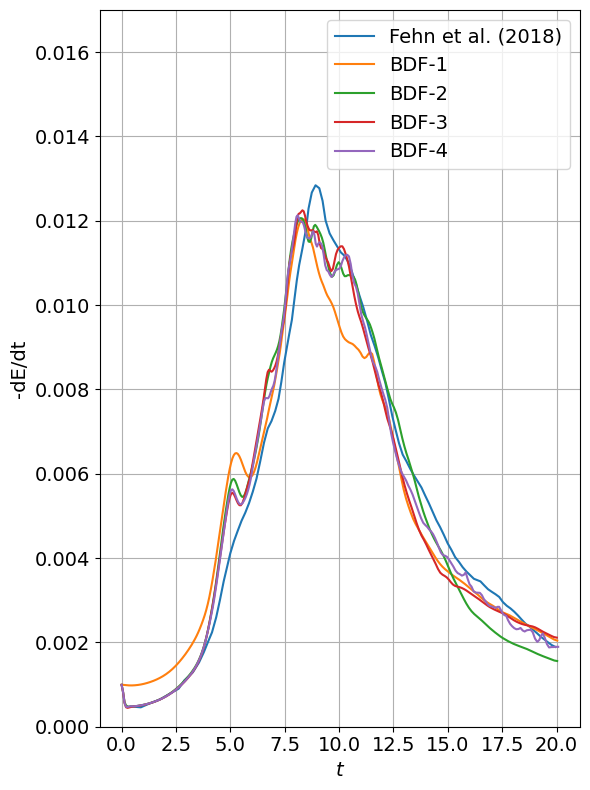}
  \caption{Kinetic energy dissipation rate}
\end{subfigure}%
\begin{subfigure}{.33\textwidth}
  \centering
  \includegraphics[width=0.99\linewidth]{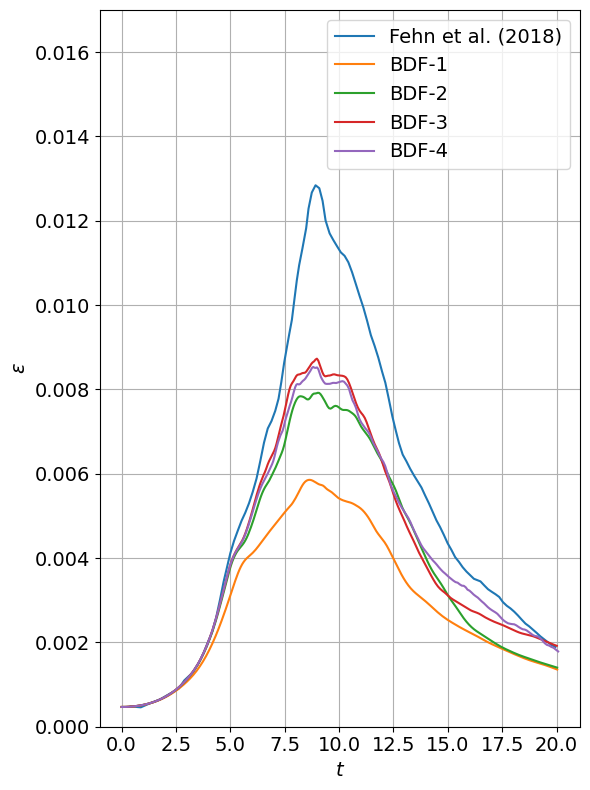}
  \caption{Molecular energy dissipation rate}
\end{subfigure}
\begin{subfigure}{.33\textwidth}
  \centering
  \includegraphics[width=0.99\linewidth]{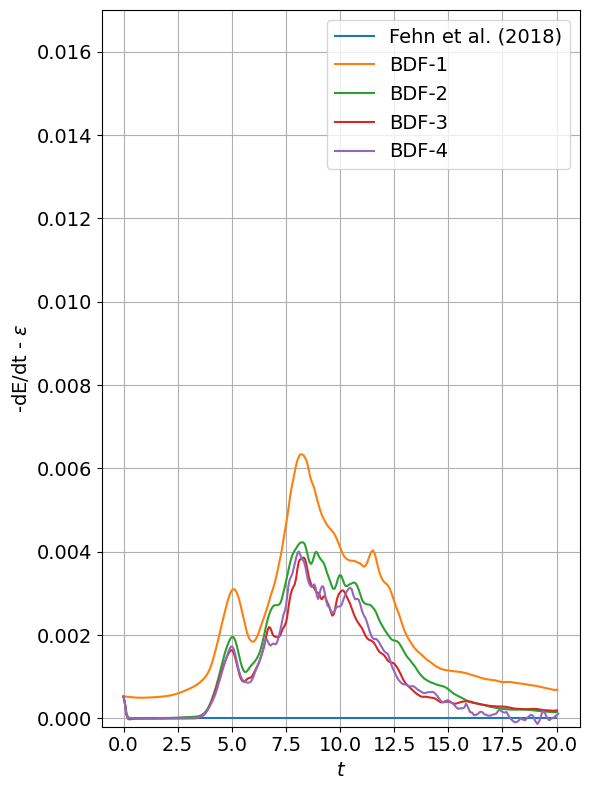}
  \caption{Numerical dissipation rate}
\end{subfigure}
    \caption{Kinetic energy, molecular dissipation and numerical dissipation of the consistent splitting method with different BDF schemes at a CFL number of 1.6. Lower numerical dissipation means more accurate, limited by the fixed spatiotemporal resolution chosen.
    }
    \label{fig:taylor_green_bdf}
\end{figure}

In order to show the effect of higher-order time integrators on the numerical dissipation, we fix the spatial discretization and target CFL number to 1.6 (stability limit for BDF-4) and increase the time integration order using the linearly implicit variant of the consistent splitting scheme. Again we choose $J_c = J_p = J$ for BDF-1 and BDF-2, while for BDF-3 and BDF-4 they are selected as $J_c = J_p = J-1$.
The results presented in Fig.~\ref{fig:taylor_green_bdf} highlight that the BDF-1 scheme exhibits higher numerical dissipation compared to the other schemes. The BDF-2 scheme is close to the higher order schemes but still shows deviations in the dissipation rate $\epsilon$. The BDF-3 and BDF-4 schemes differ for $t>14.0$ in the dissipation rate $\epsilon$, otherwise their results are rather similar, indicating that the spatial error is dominant. The increased numerical dissipation observed during the initial time steps for BDF-2 and higher-order schemes arises from the initialization procedure, which involves performing single time steps of lower-order methods until the required BDF startup data is gathered.

Finally, to compare the schemes under varying CFL numbers using the \textit{linearly implicit} variants of the convective term, we define the total numerically dissipated energy as
\begin{equation}
    \int_{t=0}^{t=20} -\frac{\mathrm{d}E}{\mathrm{d}t} - \epsilon ~\mathrm{d}t = -(E(t=20) - E(t=0))- \int_{t=0}^{t=20} \epsilon ~\mathrm{d}t
    .
\end{equation}
To assess the influence of the extrapolation order $J_c$ of the convective term, we compare the consistent splitting method for the cases $J_c = J - 1$ and $J_c = J$. For completeness, results obtained with the dual splitting scheme are also included.

As BDF-1 shows large dissipation and BDF-2 is unconditionally stable, BDF-1 is not further analyzed. Tab.~\ref{tab:taylor_green_bdf} shows the total numerically dissipated energy for BDF-2, BDF-3, and BDF-4, with the latter adopting the iterated variant of the consistent splitting scheme.

Again, the penalty terms are applied to all schemes in a postprocessing step. With the choice $J_c = J_p = J -1$, the consistent splitting method using BDF-3 remains stable at higher CFL numbers compared to the consistent splitting method and the dual splitting scheme with $J_c = J$ and $J_p = J - 1$, while producing results of similar quality.

The results indicate that selecting $J_c=J-1$ rather than $J_c = J$ in the consistent splitting method increases robustness without compromising accuracy.
A comparison between BDF-2 and BDF-3 again demonstrates the improved accuracy of the higher-order method.
For the iterated method ($\Tilde{\text{CS}}$ with $J_c = J_p = J$), an increased stability regime can be observed, compared to the iteration-free version CS$^\star$ and the dual splitting scheme. The CFL number can be increased to 2.0 for BDF-4, while the accuracy remains comparable to the iteration free version CS. Note here again that the extrapolation order for the pressure in the iteration free consistent splitting scheme is lowered, further motivating the use of the more efficient iteration-free method.
\begin{table}[ht]
  \centering
  \caption{Total numerically dissipated energy obtained via the linearly implicit consistent splitting scheme with modified pressure variable (CS, when used with $J_c = J-1$, CS$^\star$, when used with $J_c = J$ and $\Tilde{\text{CS}}$ for the iterated version) and dual splitting (DS) methods with different BDF time integration schemes.
Lower numerical dissipation means more accurate, limited by the fixed spatiotemporal resolution chosen. Simulations that fail to converge are marked by a dash.}
  \label{tab:taylor_green_bdf}
  \begin{tabular}{|l|l|c|c|c|c|}
    \hline
    \multicolumn{2}{|c|}{{CFL number}} & {0.4} & {1.0} & {2.0} & {4.0}\\
    \hline
    \multirow{2}{*}{{CS}} 
      & BDF-2 & 0.02051 & 0.02351 & 0.02873 & 0.03466\\
      & BDF-3 & 0.01983 & 0.02152 & 0.01933 & 0.01255 \\
    \hline
    \multirow{1}{*}{{CS$^\star$}} 
      & BDF-3 & 0.01985 & 0.02179 & - & -\\
       \hline
    \multirow{2}{*}{{$\Tilde{\text{CS}}$}}
      & BDF-3 & 0.02164 & 0.02465 & 0.02624 & -\\
      & BDF-4 & 0.02135 & 0.02449 & 0.01949 &  -\\
    \hline
    \multirow{2}{*}{{DS}} 
      & BDF-2 & 0.02043 & 0.02342 & 0.02846 & 0.03461\\
      & BDF-3 & 0.01982 & 0.02174 & - & - \\
    \hline
  \end{tabular}
\end{table}

In summary, the results in this section indicate that i) the scheme is higher-order accurate in space and time, ii) implicit and linearly implicit variants of the convective term allow target CFL numbers (much) larger than 1 in practical examples, iii) BDF-3 and higher can be used to improve accuracy when using larger time steps also for non-smooth solutions showcased in the 3D Taylor--Green vortex benchmark, iv) the consistent splitting scheme with $J_c = J_p = J-1$ shows improved stability for BDF-3 compared to the dual splitting scheme and similar accuracy to the coupled solution approach.

\section{Concluding Remarks}
\label{sec:conclusion}
This work proposes a discontinuous Galerkin spatial discretization of the consistent splitting scheme by~\citet{Liu2009}, with a focus on the formulation of flux terms in an $L^2$-conforming setting and the incorporation of consistent boundary conditions. The convective term is formulated in a linearly implicit manner, overcoming the strict CFL condition of an explicit formulation
while avoiding the solution of a nonlinear system, which would otherwise result from a fully implicit formulation. By choosing appropriate fluxes for the divergence of the convective operator, the forcing term and the divergence of the velocity for the PPE and integrating each contribution by part once, the consistent pressure boundary condition reduces to contributions from only the acceleration and viscous term. In comparison to~\citet{Liu2010}, the extrapolation order of the term is \textit{lowered} without impacting the temporal convergence rate as shown in section~\ref{sec:numerical_results_manufatured}, but leading to improved stability with higher order BDF schemes as observed also for the related scheme by \citet{Karniadakis1991}.

For stabilization, a continuity penalty term and a divergence penalty term are added, while Leray projection is used to improve mass conservation.
In contrast to a continuous Galerkin discretization, rewriting the curl-curl term using the Gauss theorem (on each element individually for $L^2$-conforming DG adopted herein) and replacing the second-order derivatives with first-order ones similar to~\cite{Liu2009} or \cite{Creff2025} leads to a decrease in the observed convergence rates. 

Lowest equal-order discretizations with the present scheme might therefore lead to a loss of accuracy if the curl-curl term cannot be resolved sufficiently for highly viscous problems, since the curl-curl term cannot be represented. Note here that using inf-sup stable pairs, as implicitly required by the Leray projection, being the projection onto a divergence-free space, these drawbacks seem negligible. From a performance perspective, the pressure Poisson equation often constitutes the computational bottleneck, such that a richer velocity approximation as used in an inf-sup stable approach is often preferred.

Despite these minor limitations, the proposed work offers two main advantages: i) compared to the coupled solution approach with a saddle-point linear system, velocity and pressure are decoupled, leading to smaller problems that can be effectively preconditioned by standard tools, and ii) compared to standard projection methods, splitting errors are avoided. The splitting errors often lead to spurious pressure boundary layers, introduced by artificial boundary conditions formulated to render the pressure Poisson problem solvable.  Additionally, this renders BDF-3 and BDF-4 possible choices for increased accuracy, albeit with the usual limitations due to the missing A-stability of those time integrators.

The variant considered within this work is demonstrated to yield higher-order convergence rates in space and time in an example with known analytical solution, numerically confirming the formulated claims. Further, we show excellent agreements with state-of-the-art incompressible flow solvers in practically relevant benchmarks being the flow past a cylinder problem and the Taylor--Green vortex problem, see \citet{Schaefer1996} and \citet{TaylorGreen1937, Wang2013}. In the latter example, the scheme is compared to its closest competitor, the dual splitting scheme by \citet{Karniadakis1991}, showing very much similar performance, but with the important difference being that i) the present scheme is not affected by splitting errors, and ii) extensions were recently shown to be unconditionally stable in time also for BDF-3 and even higher~\citep{Huang2025stability}---a property previously only proven for other splitting or projection schemes \textit{up to second order}, albeit at the price of potentially higher error constants. Consistent splitting schemes are thus great candidates for the numerical solution of incompressible flow problems. In this sense, the present work is a step towards the successful $L^2$-conforming discontinuous Galerkin discretization of this family of methods.

To avoid the need for additional stabilization terms, especially the continuity-stabilization, the consistent-splitting scheme can be extended to H(div) conforming elements for the velocity variables, while keeping the $L^2$ discretizations as the natural ansatz space for the PPE. This leads to questions regarding the performance of these methods taking approximation quality and time to solution into account, where ongoing developments on the preconditioning of the momentum equation in convection-dominated scenarios play a key role. 

\section*{Acknowledgements}
% PDExa
This work received support by the German Federal
Ministry of Research, Technology and Space (BMFTR) through project
``PDExa: Optimized Software Methods for Solving Partial Differential
Equations on Exascale Supercomputers'', grant agreement nos.~16ME0637K, 16ME0638, and 16ME0640
and the European Union -- NextGenerationEU,
% dealii-X
the
EuroHPC joint undertaking Centre
of Excellence dealii-X, grant agreement no.~101172493,
% DFG
by the German research foundation through project
no.~524455704 (``High-Performance Simulation Tools for Hemodynamics''),
% ERC
and by BREATHE, an ERC-2020-ADG Project, Grant Agreement ID 101021526.
% LNM support
The authors thank W.A. Wall and M. Bergbauer for the inspiring discussions and their support.

\section*{Appendix} \label{sec:Appendix}
\subsection{Similarity to the Higher-Order Dual Splitting with Periodic Boundary Conditions}
Starting from the dual splitting in time discrete formulation Equations~\eqref{eq:dual_splitting_1}--\eqref{eq:dual_splitting_helmholtz}:
\begin{align*}
  \frac{\gamma_0 \mathbf{u}^* - \sum_{i = 0}^{J-1} (\alpha_i \mathbf{u}^{n-i})}{\Delta t} &= -\sum_{i = 0}^{J-1} \beta_i (\mathbf{u}^{n-i} \cdot \nabla) \mathbf{u}^{n-i} +  \mathbf{f}^{n+1} ,\\
   -\Delta p^{n+1} &= -\frac{\gamma_0}{\Delta t} \nabla \cdot \mathbf{u}^*, \\ 
  \mathbf{u}^{**} &= \mathbf{u}^* - \frac{\Delta t}{\gamma_0} \nabla p^{n+1}, \\
 \frac{\gamma_0}{\Delta t} \mathbf{u}^{n+1} - \nu \Delta \mathbf{u}^{n+1}  &= \frac{\gamma_0}{\Delta t} \mathbf{u}^{**}.
\end{align*}
Rewriting the explicit convection step gives
\begin{align*}
      \frac{\gamma_0}{\Delta t} \mathbf{u}^*  =  -\sum_{i = 0}^{J-1} \beta_i (\mathbf{u}^{n-i} \cdot \nabla) \mathbf{u}^{n-i} +  \mathbf{f}^{n+1} +   \sum_{i = 0}^{J-1} \frac{\alpha_i}{\Delta t} \mathbf{u}^{n-i}.
\end{align*}
Inserting into the Poisson Equation~\eqref{eq:dual_splitting_pressure_poisson} leads to
\begin{align*}
     -\Delta p^{n+1} &= -\frac{\gamma_0}{\Delta t} \nabla \cdot \mathbf{u}^* = \nabla \cdot \left(\sum_{i = 0}^{J-1} \beta_i (\mathbf{u}^{n-i} \cdot \nabla) \mathbf{u}^{n-i} -  \mathbf{f}^{n+1} -   \sum_{i = 0}^{J-1} \frac{\alpha_i}{\Delta t} \mathbf{u}^{n-i}\right),
\end{align*}
which is equal to PPE~\eqref{eq:ppe_time_discretized_leray} with Leray projection.
Further, the Helmholtz like Equation~\eqref{eq:dual_splitting_helmholtz} can be written by inserting the definitions of $\ve{u}^*$ and $\ve{u}^{**}$ as
\begin{gather*}
     \frac{\gamma_0}{\Delta t} \mathbf{u}^{n+1} - \nu \Delta \mathbf{u}^{n+1}  
     = 
     \frac{\gamma_0}{\Delta t} \mathbf{u}^{**}
     =
     \frac{\gamma_0}{\Delta t}\mathbf{u}^* - \nabla p^{n+1} 
     \\
     =  -\sum_{i = 0}^{J-1} \beta_i (\mathbf{u}^{n-i} \cdot \nabla) \mathbf{u}^{n-i} +  \mathbf{f}^{n+1} +   \sum_{i = 0}^{J-1} \frac{\alpha_i}{\Delta t} \mathbf{u}^{n-i} - \nabla p^{n+1},
\end{gather*}
which gives the momentum Equation~\eqref{eq:momentum_equation} with an explicit convection term
\begin{align*}
         \frac{\gamma_0}{\Delta t} \mathbf{u}^{n+1} - \nu \Delta \mathbf{u}^{n+1} =  \mathbf{f}^{n+1}  -\sum_{i = 0}^{J-1} \beta_i (\mathbf{u}^{n-i} \cdot \nabla) \mathbf{u}^{n-i} +   \sum_{i = 0}^{J-1} \frac{\alpha_i}{\Delta t} \mathbf{u}^{n-i} - \nabla p^{n+1}.
\end{align*}

\subsection{Pressure Errors}
Table~\ref{tab:consistent_splitting_variants_pressure} lists the pressure errors and respective convergence rates complementing Table~\ref{tab:consistent_splitting_variants}.

\begin{table}[ht]
\centering
\caption{Relative $L^2$-errors and observed convergence rates in parentheses of the pressure for BDF orders 1 to 4 and various consistent splitting scheme variants as summarized in Table~\ref{tab:overview_schemes}.}
\label{tab:consistent_splitting_variants_pressure}
\resizebox{\linewidth}{!}{
\begin{tabular}{|c|c|c|c|c|c|}
\hline
Scheme & $\Delta t/T \times10^{-2}$ &  BDF-1 & BDF-2 & BDF-3 & BDF-4 \\
\hline

Leray
& $3.13$ & 1.74e-1 (1.04) &  7.06e-3 (2.30) & 2.95e-3 (2.52) & 1.94e-4 (3.39) \\
Projection & $1.65$ & 8.41e-2 (1.05) & 1.50e-3 (2.23) & 4.98e-4 (2.57) & 3.67e-5 (2.40) \\
& $0.78$ & 4.05e-2 (1.05) & 3.35e-4 (2.17) & 8.35e-5 (2.58)  & 2.41e-6 (3.93) \\
\hline
 
Modified & $3.12$ &  8.73e-2 (1.08)  & 3.62e-3 (2.10) & 1.71e-4 (3.14) & 8.81e-6 (4.14) \\
Pressure& $1.56$ &  4.22e-2 (1.05) & 8.70e-4 (2.06) & 2.04e-5 (3.07) & 5.34e-7 (4.04) \\
& $0.78$ & 2.07e-2 (1.03) & 2.13e-4 (2.03) &  2.51e-6 (3.02) & 3.93e-8 (3.76) \\
\hline

\multirow{3}{*}{Iterated} 
& $3.12$ & 4.22e-2 (1.05)  & 8.70e-4 (2.06) & 2.03e-5 (3.07) & 5.03e-7 (4.10) \\
& $1.56$ & 2.07e-2 (1.03) & 2.13e-4 (2.03) & 2.47e-6 (3.04)  & 2.98e-8 (4.08) \\
& $0.78$ & 1.02e-2 (1.02) & 5.26e-5 (2.02) & 3.05e-7 (3.02)  & 5.95e-9 (2.33) \\
\hline

\multirow{3}{*}{Traction} 
& $3.12$ & 7.63e-2 (1.23) & 2.43e-3 (3.70) & 8.76e-4 (4.65) & 8.22e-5 (6.08) \\
& $1.56$ & 3.48e-2 (1.13) & 4.99e-4 (2.28) & 1.09e-4 (3.01) & 4.59e-6 (4.16) \\
& $0.78$ & 1.65e-2 (1.07) & 1.11e-4 (2.16) & 1.36e-5 (3.01) & 3.85e-7 (3.58) \\
\hline

\end{tabular}}
\end{table}

\bibliographystyle{unsrtnat} 
\bibliography{references}

\end{document}